\documentclass[MSNbibl,number,citesort,seceqn,dvips]{arxbj}
\usepackage{dcolumn}

\aid{0}
\volume{20}
\issue{3}
\pubyear{2014}
\firstpage{1344}
\lastpage{1371}
\doi{10.3150/13-BEJ524} %kopijuoti is PTS

\makeatletter

\newcolumntype{d}[1]{D{.}{.}{#1}}
\newcommand{\rrvert}{\vert}
\newcommand{\llvert}{\vert}
\newcommand{\eqref}[1]{(\ref{#1})}
\def\I{\mathbf{1}}
\def\m{\scriptsize \maltese}
\def\d{\mathrm{d}}

\newcommand{\com}[1]{#1}

\newtheorem{thmm}{Theorem}[section]
\newtheorem{prop}{Proposition}[section]
\newproclaim{defn}{Definition}[section]
\newtheorem{lem}{Lemma}[section]
\newtheorem{cor}{Corollary}[section]
\newremark{rem}{Remark}[section]
\newproclaim{algo}{Algorithm}[section]
\makeatother

\begin{document}
\begin{frontmatter}

\title{On the empirical multilinear copula process for count data}
\runtitle{Empirical copula process for count data}

\begin{aug}
%%%% inicialai - be tarpu
\author[a]{\inits{C.}\fnms{Christian} \snm{Genest}\thanksref{a,e1,u1}\ead[label=e1,mark]{cgenest@math.mcgill.ca}\ead[label=u1,url,mark]{www.math.mcgill.ca/cgenest/}},
\author[a]{\inits{J.G.}\fnms{Johanna G.} \snm{Ne{\v s}lehov\'a}\ead[label=e2,mark]{johanna@math.mcgill.ca}\ead[label=u2,url,mark]{www.math.mcgill.ca/neslehova/}\thanksref{a,e2,u2}}
\and\\
\author[b]{\inits{B.}\fnms{Bruno} \snm{R\'emillard}\thanksref{b}\ead[label=e3]{bruno.remillard@hec.ca}\ead[label=u3,url]{neumann.hec.ca/pages/bruno.remillard/}}
\runauthor{C. Genest, J.G. Ne{\v s}lehov\'a and B. R\'emillard} %% auto
\address[a]{%Christian Genest \& Johanna G. Ne{\v s}lehov\'a\\
Department of Mathematics and Statistics,
McGill University,
805, rue Sherbrooke ouest,
Montr\'eal (Qu\'ebec),
Canada H3A 0B9.\\
\printead{e1,e2};\\
\printead{u1,u2}}
\address[b]{%Bruno R\'emillard\\
Service de l'enseignement des
m\'ethodes quantitatives de gestion,
HEC Montr\'eal,
3000, chemin de la C\^{o}te-Sainte-Catherine,
Montr\'eal (Qu\'ebec),
Canada H3T 2A7.\\
\printead{e3,u3}}
\end{aug}

% HISTORY:
\received{\smonth{7} \syear{2012}}
\revised{\smonth{2} \syear{2013}}

% ABSTRACT
%
\begin{abstract}
Continuation refers to the operation by which the cumulative
distribution function of a discontinuous random vector is made
continuous through multilinear interpolation. The copula that results
from the application of this technique to the classical empirical
copula is either called the multilinear or the checkerboard copula. As
shown by Genest and
Ne{\v{s}}lehov\'a (\textit{Astin Bull.} \textbf{37} (2007) 475--515) and Ne{\v{s}}lehov{\'a}
(\textit{J. Multivariate Anal.} \textbf{98} (2007) 544--567), this
copula plays a central role in characterizing dependence concepts in
discrete random vectors. In this paper, the authors establish the
asymptotic behavior of the empirical process associated with the
multilinear copula based on $d$-variate count data. This empirical
process does not generally converge in law on the space $\mathcal
{C}([0,1]^d)$ of continuous functions on $[0,1]^d$, equipped with the
uniform norm. However, the authors show that the process converges in
$\mathcal{C}(K)$ for any compact $K \subset\mathcal{O}$, where
$\mathcal{O}$ is a dense open subset of $[0,1]^d$, whose complement is
the Cartesian product of the ranges of the marginal distribution
functions. This result is sufficient to deduce the weak limit of many
functionals of the process, including classical statistics for monotone
trend. It also leads to a powerful and consistent test of independence
which is applicable even to sparse contingency tables whose dimension
is sample size dependent.
\end{abstract}

% KEYWORDS
% visi is mazosios raides ir pagal abecele
%

\end{frontmatter}

%s1 #&#
\section{Introduction}\label{sec1}
This paper's central message is that there are advantages, both
conceptual and technical, to viewing a contingency table as arising
from a multivariate distribution having uniform margins on the unit
interval, that is, a copula. As will be shown here, this approach leads
to new statistical methodology that can be used to analyze tables that
are sparse or whose number of categories grows with the sample size.

To go straight to the point, consider the simple case of a $K \times L$
contingency table derived from a random sample of size $n$ of ordinal
or interval responses in ordered categories $A_1 < \cdots< A_K$ and
$B_1 < \cdots< B_L$. For arbitrary $k \in\{1, \ldots, K\}$ and $\ell
\in\{1, \ldots, L\}$, let $f_{k\ell}$ be the relative frequency of
the pair $(A_k, B_\ell)$ and denote by $f_{k +}$ and $f_{+ \ell}$ the
row-wise and column-wise totals, respectively. Further set
\[
F_{k+} = \sum_{i=1}^k
f_{i+},\qquad F_{+\ell} = \sum_{i=1}^\ell
f_{+i},
\]
and let $F_{0+} = F_{+0} = 0$. A density $\hat c_n^{\m}$ with respect to
the Lebesgue measure can then be defined (almost everywhere) on
$[0,1]^2$ by setting
\[
\hat c_n^{\m}(u,v) = \frac{f_{k\ell}}{f_{k+} f_{+\ell}}
\]
whenever\vspace*{1pt} $u \in(F_{(k-1)+}, F_{k+})$ and $v \in(F_{+(\ell-1)},
F_{+\ell})$. As shown in Section~\ref{sec2}, the corresponding
distribution function $\widehat C_n^{\m}$ is a copula, that is, its
margins are uniform on $[0,1]$. Moreover, when $f_{k\ell} =
f_{k+}f_{+\ell}$ for all $k\in\{1,\ldots, K\}$ and $\ell\in\{
1,\ldots, L\}
$, $\widehat C_n^{\m}$ becomes the independence copula $\Pi$ whose
Lebesgue density is identically equal to $1$ on $[0,1]^2$.

More significantly, several standard measures of association in the
pair $(X,Y)$, and classical tests of independence between $X$ and $Y$,
are based on $\widehat C_n^{\m}$. For example, Pearson's $\chi^2$
statistic and the likelihood ratio statistic $G^2$ are immediately seen
to satisfy
%
%
%e1.1 #&#
\begin{eqnarray}
\chi^2 &=& n \sum_{k=1}^K \sum
_{\ell=1}^L \frac{ (f_{k\ell} -
f_{k+}f_{\ell+})^2}{f_{k+}f_{\ell+}} = n\int
_0^1 \int_0^1
\bigl\{ \hat c_n^{\m} (u,v) -1 \bigr\}^2 \,\d v
\,\d u, \label{eq11}
\\
G^2 &=& 2 n \sum_{k=1}^K \sum
_{\ell=1}^L f_{k\ell} \ln\biggl(
\frac
{f_{k\ell}}{f_{k+} f_{+\ell}} \biggr) = 2n \int_0^1 \int
_0^1 \ln\bigl\{ \hat c_n^{\m}(u,v)
\bigr\} \,\d\widehat C_n^{\m}(u,v). \label{eq12}
\end{eqnarray}
With some additional work (Ne{\v{s}}lehov{\'a} \cite{Neslehova2007})
it can also be shown
that the well-known Spearman and Kendall statistics for testing
monotone trend (Agresti \cite{Agresti2007}) can be rewritten in terms of
$\widehat C_n^{\m}$. Many other examples could be given. %These and
%other
%examples are further explored in Section \ref{sec5} in the more
%general case of multi-way tables.

The introduction of the \emph{multilinear empirical copula} $\widehat
C_n^{\scriptsize\maltese}$ in this context is not merely a neat way of unifying
various known statistics for frequency data analysis. Because integral
expressions such as $\eqref{eq11}$ and $\eqref{eq12}$ make sense
even when the number of categories changes with $n$, $\widehat C^{\m}_n$
is rather a key tool for the investigation of new or existing
procedures that can be used even in cases where the table is sparse or
of varying dimension.

Further, it may be seen that when $X$ and $Y$ are continuous, $\widehat
C_n^{\scriptsize\maltese}$ is a smoothed version of the classical empirical copula
(Deheuvels \cite{Deheuvels1979}) from which it differs by at most \com
{a factor
of} $1/n$ uniformly. Statistical tools based on $\widehat C_n^{\scriptsize\maltese}$
can thus bridge the gap between continuous and discrete outcomes. In
particular, the problems associated with ties, which invalidate many of
the procedures developed for continuous data (Genest, Ne\v{s}lehov\'a and
Ruppert~\cite
{GenestNeslehovaRuppert2011}), are then automatically taken care of.
\com{While it seems intuitively reasonable to base inference on
$\widehat C_n^{\scriptsize\maltese}$, this new approach generally requires the
knowledge of its limit $C^{\scriptsize\maltese}$ and the asymptotic behavior of the
corresponding empirical process
%
%
%e1.3 #&#
\begin{equation}
\label{eq13} \mathbb{\widehat C}^{\m}_n = \sqrt{n} \bigl(
\widehat C_n^{\m}- C^{\m} \bigr),
\end{equation}
which has hitherto never been studied in the literature.}

This paper contributes to the problem by determining the asymptotic
behavior of the process \eqref{eq13} in general dimension $d \ge2$
when the components of the underlying random vector $\mathbf
{X}=(X_1, \ldots, X_d)$ are either integer-valued or strictly
increasing transformations thereof. As will be seen, $\widehat C_n^{\m}$
is a consistent estimator of the so-called multilinear extension (or
checkerboard) copula $C^{\m}$ of $\mathbf{X}$. This limiting copula,
defined in Section~\ref{sec2}, has been studied earlier, for example,
by Genest and
Ne{\v{s}}lehov\'a \cite{GenestNeslehova2007} and Ne{\v{s}}lehov{\'a}
\cite{Neslehova2007}, who showed
that it captures many important dependence properties of $\mathbf
{X}$ when $d=2$. In particular, when the components of $\mathbf{X}$
are independent, $C^{\m}$ is the independence copula $\Pi$.

The main result, stated in Section~\ref{sec3}, gives the asymptotic
behavior of the process \eqref{eq13}.
Unless the components of $\mathbf{X}$ are mutually independent,
$\mathbb{\widehat C}^{\m}_n$ does not generally converge on the space
$\mathcal{C}([0,1]^d)$ of continuous functions on $[0,1]^d$ equipped
with the uniform norm because $C^{\m}$ has discontinuous partial
derivatives. Fortunately, $\mathbb{\widehat C}^{\m}_n$
converges~-- without any regularity conditions~-- in the subspace
$\mathcal{C}(K)$ for any compact subset $K\subset\mathcal{O}$, where
\com{$\mathcal{O}$} is a dense open subset of $[0,1]^d$ \com{whose
complement is the Cartesian product of the ranges of the marginal
distribution functions}. The proof of the main result is involved; it
is outlined in Section~\ref{sec4} and detailed in the \hyperref[app]{Appendix}.

\com{To illustrate the usefulness of the process \eqref{eq13} for
inference, Section~\ref{sec5} provides a few initial examples of
application. It is first shown that the main result is sufficient to
deduce the limiting distribution of classical statistics for monotone
trend such as Spearman's rho and Kendall's tau. Moreover, a new and
consistent Cram\'er--von Mises type test of independence is proposed
that can be used whatever the margins. As illustrated through a small
simulation study, it performs very well even for sparse contingency
tables whose dimension is sample size dependent; in all cases
considered, it is consistently more powerful than the classical
chi-squared test.} Section~\ref{sec6} concludes.

%s2 #&#
\section{The multilinear extension copula} \label{sec2}

Suppose that $\mathbf{X}=(X_1,\ldots,X_d)$ is a vector of discrete
random variables with joint cumulative distribution function $H$. For
each $j \in\{ 1, \ldots, d\}$, let $F_j$ denote the distribution
function of $X_j$ and assume that \com{there exists a strictly
increasing function $A_j\dvt\mathbb{N} \to\mathbb{R}$ such that
$\operatorname
{supp}(X_j) \subseteq\{ A_j (k)\dvt k \in\mathbb{N} \}$. Note that the
inclusion may be strict; in particular, it is not assumed that $\Pr\{
X_j = A_j(k) \} > 0$ holds for all $k \in\mathbb{N}$ or that the
support of $X_j$ is infinite. Furthermore, observe that the closure of
the range of $F_j$, viz. $\mathcal{R}_j = \{ 0,1, F_j\{ A_j (0)\},
F_j\{ A_j (1)\} , \ldots\}$, defines a partition of $[0,1]$. In what
follows, $A_j (-1) = A_j(0) - 1$ for all $j \in\{ 1, \ldots, d \}$ by
convention.}

%
%de2.1 #&#
\begin{defn}
The multilinear extension copula $C^{\scriptsize\maltese}$ of $H$ is the unique
copula whose density \com{with respect to the Lebesgue measure} is
given by
\[
c^{\m}(u_1,\ldots, u_d) = \frac{\Pr\{X_1=A_1(k_1),\ldots, X_d =
A_d(k_d)\}
}{\Pr\{X_1=A_1(k_1)\}\times\cdots\times\Pr\{X_d=A_d(k_d)\}}
\]
whenever for all $j\in\{1,\ldots, d\}$, $F_j\{A_j(k_j-1)\} < u_j \le
F_j\{A_j(k_j)\}$ for some $k_j \in\mathbb{N}$.
\end{defn}

An explicit form of $C^{\m}$, which is easily verified by
differentiation, is given in Proposition \ref{prop21} below. For each
$j \in\{ 1, \ldots, d\}$ and $u \in[0,1]$, let $u_j^-$ and $u_j^+$
be, respectively, the greatest and the least element of $\mathcal{R}_j$
such that $u_j^- \le u \le u_j^+$. Further let
\[
\lambda_{F_j}(u) = %
\cases{ \bigl(u - u_j^-
\bigr)/ \bigl(u_j^+ - u_j^- \bigr), & \quad$\mbox{if }
u_j^- \neq u_j^+$, \vspace*{2pt}
\cr
1 , & \quad $\mbox{otherwise}.$}
\]
Thus when $k \in\mathbb{N}$ is such that $\Delta F_j\{A_j(k)\} = \Pr
\{
X_j = A_j(k)\} > 0$,
then for all $u \in(F_j\{A_j(k-1)\}, F_j\{A_j(k)\})$, one has
$u_j^-=F_j\{A_j(k-1)\}$, $u_j^+ = F_j\{A_j(k)\}$ and
\[
\lambda_{F_j}(u) = \frac{u-F_j\{A_j(k-1)\}}{\Delta F_j\{A_j(k)\}}.
\]
Furthermore, if $F_j^{-1}$ is the pseudo-inverse of $F_j$, then $F_j
\circ F_j^{-1}(u_j^-) = F_j \{A_j(k-1)\}$ and $F_j \circ
F_j^{-1}(u_j^+) = F_j \{A_j(k)\}$. Finally, for
any $S \subset\{ 1, \ldots, d\}$ and $u_1 ,\ldots, u_d \in[0,1]$, set
\[
\lambda_{H,S} (u_1, \ldots, u_d) = \prod
_{\ell\in S}\lambda_{F_\ell
}(u_\ell) \prod
_{\ell\notin S} \bigl\{1-\lambda_{F_\ell}(u_\ell)
\bigr\},
\]
\com{which depends on $H$ only through its margins $F_1,\ldots, F_d$.}
%
%
%pr2.1 #&#
\begin{prop}
\label{prop21}
The multilinear extension copula $C^{\m}$ of $H$ is given by
\[
C^{\m}(u_1,\ldots, u_d) = \sum
_{S\subseteq\{1,\ldots,d\}} \lambda_{H,S} (u_1, \ldots,
u_d) H \bigl\{F_1^{-1}(u_{S_1}),\ldots,
F^{-1}_d(u_{S_d}) \bigr\},
\]
where for each $j \in\{1,\ldots, d\}$, $u_{S_j} = u_j^+$ if $j \in S$
and $u_{S_j} = u_j^-$ otherwise. \com{In particular, $C^{\m}(u_{S_1},
\ldots, u_{S_d}) = H \{F_1^{-1}(u_{S_1}),\ldots,
F^{-1}_d(u_{S_d}) \}$ for any $S\subseteq\{1,\ldots,d\}$.
}
\end{prop}

It is easily seen that $C^{\m}$ satisfies Sklar's representation, that
is, for all $x_1,\ldots, x_d \in\mathbb{R}$,
\[
H(x_1, \ldots, x_d) = C^{\m} \bigl
\{F_1(x_1), \ldots, F_d(x_d) \bigr
\}.
\]
This is because in effect, this identity needs only be verified if for
all $j \in\{ 1, \ldots, d\}$, $x_j = A_j(k_j)$ for some $k_j \in
\mathbb{N}$ such that $\Delta F_j\{A_j(k_j)\} >0$. In fact, $C^{\m}$ is
precisely the construction used to extend a sub-copula to a copula in
the proof of Sklar's theorem; see, for example, Nelsen \cite
{Nelsen1999} for
details in the bivariate case.

The copula $C^{\m}$ is known to capture many important dependence
properties of $H$, as summarized by Genest and
Ne{\v{s}}lehov\'a \cite{GenestNeslehova2007}. As
shown by Ne{\v{s}}lehov{\'a} (\cite{Neslehova2007}, Corollary 6),
$C^{\m}$ is invariant with
respect to strictly increasing transformations of the margins.

Now consider a random sample $\mathcal{X} = \{ (X_{11}, \ldots, X_{1d}),
\ldots, (X_{n1}, \ldots, X_{nd}) \}$ from $H$ and let $H_n$ be the
corresponding empirical distribution function. \com{Because $H_n$ is
itself a discrete distribution, one can define its multilinear
extension copula $\widehat C_n^{\m}$ and its corresponding density
$\hat
c^{\m}_n$ with respect to the Lebesgue measure as above. To be explicit,
fix $j \in\{1,\ldots, d\}$ and denote by $A_{nj}(0) < \cdots<
A_{nj}(n_j)$ the distinct values of $X_{1j},\ldots, X_{nj}$. Let also
$A_{nj}(-1) = A_{nj}(0)-1$. The range $\mathcal{R}_{nj}$ of $F_{nj}$
then consists of
\[
0=F_{nj} \bigl\{A_{nj}(-1) \bigr\} < F_{nj} \bigl
\{A_{nj}(0) \bigr\} < \cdots< 1=F_{nj} \bigl\{
A_{nj}(n_j) \bigr\}.
\]
If $(u_1,\ldots, u_d)$ is such that for all $j\in\{1,\ldots, d\}$,
$F_{nj}\{A_{nj}(k_j-1)\} < u_j \le F_{nj}\{A_{nj}(k_j)\}$ for some $k_j
\in\{0,\ldots, n_j\}$, then
\[
\hat c_n^{\m}(u_1,\ldots, u_d) =
\frac{h_n\{A_{n1}(k_1),\ldots,A_{nd}(k_d)\}
}{\Delta F_{n1} \{A_{n1}(k_1)\}\times\cdots\times\Delta F_{nd} \{
A_{nd}(k_d)\}} ,
\]
whose numerator is the proportion of data with $X_{ij} = A_{nj}(k_{j})$
for $j\in\{1,\ldots, d\}$, and}
%
%
%e2.1 #&#
\[
\label{eq21} \widehat C_n^{\m}(u_1,\ldots,
u_d)
= \sum_{S\subseteq\{1,\ldots,d\}} \lambda_{H_n,S}
(u_1, \ldots, u_d) H_n \bigl
\{F_{n1}^{-1}(u_{S_1}),\ldots, F^{-1}_{nd}(u_{S_d})
\bigr\} .
\]

Observe that\vspace*{1pt} $\hat c_n^{\m}$ and $\widehat C_n^{\m}$ are both
functions of
the component-wise ranks. As announced in the \hyperref
[sec1]{Introduction}, $\widehat
C_n^{\m}$ is a consistent estimator of the multilinear extension copula
$C^{\m}$ of $H$. This fact will be a consequence of this paper's main
result, Theorem \ref{thm31}, which characterizes the limit of the
process $\mathbb{\widehat C}^{\m}_n$ defined in \eqref{eq13}.

%
%re2.1 #&#
\begin{rem}
\label{rem21}
When $X_1,\ldots, X_d$ are continuous, $\widehat C_n^{\m}$ was actually
used by Deheuvels \cite{Deheuvels-Raoult1979} to construct tests of
independence. It is then asymptotically equivalent to the empirical
copula $\widehat C_n$ given, for all $u_1,\ldots, u_d \in[0,1]$, by
\[
\widehat C_n(u_1,\ldots, u_d) =
\frac{1}{n} \sum_{i=1}^n \I\bigl\{
F_{n1}(X_{i1}) \le u_1,\ldots,F_{nd}(X_{id})
\le u_d \bigr\}.
\]
Indeed, if for all $j\in\{1,\ldots, d\}$, \com{$F_{nj}\{
A_{nj}(k_j-1)\}
\le u_j < F_{nj}\{A_{nj}(k_j)\}$} for some \com{$k_j \in\{0,\ldots
,n_j\}$}, then \com{$\widehat C_n(u_1,\ldots, u_d) =H_n\{
A_{n1}(k_1-1),\ldots,A_{nd}(k_d-1)\}$}. Because the coefficients
$\lambda
_{H_n,S}$ are non-negative and add up to $1$ by the multinomial formula,
the fact that $H_n$ is non-decreasing component-wise implies that
\[
\com{\widehat C_n(u_1,\ldots, u_d) \le
\widehat C_n^{\m}(u_1,\ldots, u_d)
\le H_n \bigl\{A_{n1}(k_1),\ldots,A_{nd}(k_d)
\bigr\}}.
\]
Hence, $| \widehat C_n(u_1,\ldots, u_d) - \widehat C_n^{\m}(u_1,\ldots,
u_d)|$ is bounded above by
\begin{eqnarray*}
&& H_n \bigl\{A_{n1}(k_1),\ldots,A_{nd}(k_d)
\bigr\} - H_n \bigl\{A_{n1}(k_1-1),\ldots
,A_{nd}(k_d-1) \bigr\}
\\
&&\quad \le\sum_{j=1}^d \bigl| F_{nj} \bigl
\{A_{nj}(k_j) \bigr\} - F_{nj} \bigl\{
A_{nj}(k_j-1) \bigr\}\bigr|,
\end{eqnarray*}
from which it follows that $\|\widehat C_n - \widehat C_n^{\m}\| \le d/n$
almost surely. \com{This also implies that $\widehat C_n^{\scriptsize\maltese}$ is
asymptotically equivalent to other versions of the empirical copula
commonly used in the literature; see, for example, Fermanian,
Radulovi{\'c} and
Wegkamp \cite
{FermanianRadulovicWegkamp2004}.}
\end{rem}

To ease the notation, it will be assumed henceforth, without loss of
generality, that $X_1, \ldots, X_d$ are integer-valued. In this case,
one has the following alternative representation of $C^{\m}$, which is
useful to study the process \eqref{eq13}.
%
%
%pr2.2 #&#
\begin{prop}
\label{prop22}
Let $(X_1,\ldots, X_d)$ be a random vector in $\mathbb{N}^d$ with
distribution function~$H$. Let also $U_1,\ldots, U_d$ be independent
standard uniform random variables, independent of $(X_1, \ldots, X_d)$.
Then $C^{\m}$ is the unique copula of the distribution function $H^{\m
}$ of
$(X_1 + U_1 -1, \ldots, X_d + U_d -1)$ with margins $F_1^{\m},\ldots,
F_d^{\m}$, that is, for all $u_1,\ldots, u_d \in[0,1]$,
\[
C^{\m}(u_1,\ldots, u_d) = H^{\m}
\bigl\{ F_1^{\m-1}(u_1),\ldots,
F_d^{\m-1}(u_d) \bigr\}.
\]
\end{prop}

Given an empirical distribution function $H_n$ based on a random sample
from a multivariate integer-valued distribution $H$, one can proceed as
in Proposition \ref{prop22} to define a multilinear extension $H_n^{\m}
$ whose margins $F_{n1}^{\m}, \ldots, F_{nd}^{\m}$ are continuous
extensions of the margins $F_{n1}, \ldots, F_{nd}$ of $H_n$. Furthermore,
\[
\widehat C_n^{\m}(u_1,\ldots, u_d)
= H_n^{\m} \bigl\{ F_{n1}^{\m-1}(u_1),
\ldots, F_{nd}^{\m-1}(u_d) \bigr\}
\]
holds for all $u_1, \ldots, u_d \in[0,1]$, which will come in handy
in Section~\ref{sec3}.

%s3 #&#
\section{The empirical multilinear copula process}
\label{sec3}

In what follows, $\mathcal{C}(K)$ stands for the space of all
continuous functions from a compact set $K \subseteq[0,1]^d$ to
$\mathbb{R}$ equipped with the uniform norm, that is, $\| f \|_K =
\sup
\{ |f(u_1,\ldots, u_d)|\dvt(u_1,\ldots,  u_d) \in K\} $. When $K = [0,1]^d$,
the index on $\| \cdot\|$ is suppressed. Similarly, let $\ell^\infty
(K)$ denote the space of all bounded functions from $K$ to $\mathbb{R}$
equipped with the uniform norm. For each $j \in\{1, \ldots, d \}$ and
all $u_1, \ldots, u_d \in(0,1)$ \com{where the partial derivatives
exist, set}
\[
\dot C_j^{\m}(u_1, \ldots, u_d) =
\frac{\partial}{\partial u_j} C^{\m} (u_1, \ldots, u_d).
\]
Furthermore, let $\mathbb{B}_{C^{\m}}$ be a $C^{\scriptsize\maltese}$-Brownian bridge,
that is, a centred  Gaussian process on $[0,1]^d$ with covariance given,
for all
$s_1, \ldots, s_d, t_1, \ldots, t_d \in[0,1]$, by
\[
C^{\m}(s_1 \wedge t_1, \ldots, s_d
\wedge t_d) - C^{\m}(s_1, \ldots,
s_d) C^{\m}(t_1, \ldots, t_d).
\]
Here, $a \wedge b = \min(a,b)$ for arbitrary $a, b \in\mathbb{R}$.
The limit of $\mathbb{\widehat C}^{\m}_n$ can be expressed in terms
of a
transformation of $\mathbb{B}_{C^{\m}}$ involving the following operator.

%
%de3.1 #&#
\begin{defn}
\label{def31}
Let $H$ be a multivariate distribution function with support included
in $\mathbb{N}^d$ and margins $F_1,\ldots, F_d$. The multilinear
interpolation operator $\mathfrak{M}_H\dvt\ell^\infty([0,1]^d) \to
\ell
^\infty([0,1]^d) \dvtx  g \mapsto\mathfrak{M}_H(g)$ is defined,
for every $ g \in\ell^\infty([0,1]^d)$, by
\[
\mathfrak{M}_H (g) (u_1,\ldots, u_d) = \sum
_{S \subseteq\{1,\ldots,d\}} \lambda_{H,S} (u_1,
\ldots, u_d) g (u_{S_1}, \ldots, u_{S_ d}).
\]
\end{defn}
\com{As was the case with $\lambda_{H,S}$, the operator $\mathfrak
{M}_H$ depends on $H$ only through its margins.}
Although the paths of the process $\mathbb{\widehat C}_n^{\m}$ are
continuous on $[0,1]^d$ for every $n$, it cannot possibly converge in
$\mathcal{C}([0,1]^d)$ in general. This is because unless $C^{\m}=
\Pi$,
its partial derivatives exist only on
the open set
\[
\mathcal{O} = \bigcup_{(k_1,\ldots,k_d) \in\mathbb{N}^d} \bigl(F_1(k_1
-1), F_1(k_1) \bigr) \times\cdots\times
\bigl(F_d(k_d -1), F_d(k_d)
\bigr).
\]
Fortunately, the convergence of $\mathbb{\widehat C}_n^{\m}$ can be
established in $\mathcal{C}(K)$ for any compact $K \subset\mathcal
{O}$. The symbol $\rightsquigarrow$ is used henceforth to denote weak
convergence.

%
%th3.1 #&#
\begin{thmm}
\label{thm31}
Let $\mathbb{C}^{\m}= \mathfrak{M}_H(\mathbb{B}_{C^{\m}})$ and
\com{let}
$K$ be any compact subset of $\mathcal{O}$. Then, as $n\to\infty$,
$\mathbb{\widehat C}^{\m}_n \rightsquigarrow\mathbb{\widehat C}^{\m
}$ in
$\mathcal{C}(K)$, where, for all $(u_1, \ldots, u_d) \in\mathcal{O}$,
\[
\mathbb{\widehat C}^{\m}(u_1, \ldots, u_d) =
\mathbb{C}^{\m}(u_1,\ldots, u_d) - \sum
_{j=1}^d {\dot C}_j^{\m}(u_1,
\ldots, u_d) \mathbb{C}^{\m} (1, \ldots, 1,
u_j, 1, \ldots,1).
\]
\end{thmm}

This theorem can be strengthened when $X_1,\ldots,X_d$ are mutually
independent, which is the case if and only if $C^{\scriptsize\maltese}$ is the
independence copula $\Pi$.
%
%
%co3.1 #&#
\begin{cor}\label{cor31} Suppose that $C^{\m}= \Pi$. Then, as $n \to
\infty$, $\mathbb{\widehat C}^{\m}_n \rightsquigarrow\mathbb
{\widehat
C}^{\m}$ in $\mathcal{C}([0,1]^d)$.
\end{cor}
%
%
%re3.1 #&#
\begin{rem}
\label{rem31} When $X_1,\ldots, X_d$ are continuous, $C^{\m}= C$ is the
unique copula of $H$ and $\widehat C_n^{\m}$ is asymptotically equivalent
to the empirical copula $\widehat C_n$ by Remark \ref{rem21}.\vspace*{1pt} R\"
uschendorf \cite
{Rueschendorf1976} showed that under suitable regularity conditions on
$C$, $\mathbb{\widehat C}^{\m}_n \rightsquigarrow\mathbb{\widehat
C}$ as
$n \to\infty$, where $\mathbb{\widehat C}$ is defined in terms of a
$C$-Brownian bridge $\mathbb{B}_C$, for all $u_1, \ldots, u_d \in
[0,1]$, by
\[
\mathbb{\widehat C}(u_1, \ldots, u_d)=
\mathbb{B}_C (u_1, \ldots, u_d) - \sum
_{j=1}^d \dot C_j(u_1,
\ldots, u_d) \mathbb{B}_{C}(1, \ldots, 1,
u_j, 1, \ldots,1).
\]
This result has since been refined in various ways; see Segers \cite
{Segers2012} and references therein.
\end{rem}

%s4 #&#
\section{Proof of Theorem \texorpdfstring{\protect\ref{thm31}}{3.1}}
\label{sec4}

The proof of the main result is quite involved. It rests on a series of
steps and propositions that are described below. All proofs may be
found in Appendix~\ref{appB}.

Because $C^{\m}$ is a copula of $H$, it can be assumed without loss of
generality that the sample $\mathcal{X}$ from $H$ arises from a random
sample $\mathcal{V} = \{ (V_{11}, \ldots, V_{1d}), \ldots,(V_{n1},
\ldots,V_{nd}) \}$ from $C^{\m}$, that is, for every $i \in\{ 1,
\ldots, n \}$, one has $X_{i1} = F_1^{-1}(V_{i1}), \ldots, X_{id} =
F_d^{-1}(V_{id})$. If $B_n$ denotes the empirical distribution function
of this latent sample $\mathcal{V}$, it is well known that as $n \to
\infty$, the corresponding empirical process $\mathbb{B}_n = \sqrt{n}
(B_n - C^{\m})$ converges weakly in $\ell^\infty([0,1]^d)$ to the
$C^{\m}
$-Brownian bridge $\mathbb{B}_{C^{\m}}$ (van~der Vaart and
Wellner \cite{VaartWellner1996}).

The first step consists of considering the case where the margins of
$H$ are known. In contrast to the continuous case, the variables
$F_1(X_1),\ldots, F_d(X_d)$ are not uniform and their joint distribution
function $D$ is not a copula. Observe that $C^{\m}= \mathfrak{M}_H(D)$
and introduce $C_n^{\m}= \mathfrak{M}_H(D_n)$, where $D_n$ denotes the
empirical distribution function of the transformed data $ (F_1(X_{11}),
\ldots, F_d(X_{1d})), \ldots, (F_1(X_{n1}), \ldots, F_d(X_{nd}))$.
\com{Note that $C_n^{\m}$ cannot be computed in practice, because it
relies on the unknown marginal distribution functions.}
As is easily seen by differentiation, $C_n^{\m}$ is a continuous
distribution function on $[0,1]^d$ whose $j$th margin is given, for all
$u\in[0,1]$, by
\[
C_{nj}^{\m}(u) = \lambda_{F_j}(u) D_{nj}
\bigl(u^+ \bigr) + \bigl\{1-\lambda_{F_j}(u) \bigr\} D_{nj}
\bigl(u^- \bigr).
\]
Because its margins are not uniform, $C_n^{\m}$ is not a copula. The
following proposition shows that the empirical process $\mathbb
{C}_n^{\m}
= \sqrt{n} (C_n^{\m}-C^{\m})$ converges. Its proof rests on the fact
that $\mathfrak{M}_H$ is a continuous linear contraction. This is
because the weights $\lambda_{H,S}$ are non-negative and add up to~$1$,
so that for any $g, g^* \in\ell^\infty([0,1]^d)$, one has $\|
\mathfrak
{M}_H (g) - \mathfrak{M}_H(g^*) \| \le\| g - g^*\|$.

%
%pr4.1 #&#
\begin{prop}\label{prop41}
As $n \to\infty$, $ \mathbb{C}^{\m}_n \rightsquigarrow\mathbb
{C}^{\m}
=\mathfrak{M}_H(\mathbb{B}_{C^{\m}})$ in $\mathcal{C}([0,1]^d)$.
\end{prop}

Next, the process $\mathbb{\widehat C}_n^{\m}$ in which margins are
unknown can be written in the form
%
%
%e4.1 #&#
\begin{equation}
\label{eq41} \mathbb{\widehat C}_n^{\m}= \mathbb{
\widetilde C}_n^{\m} + \mathbb{\widetilde D}_n,
\end{equation}
where the summands are defined, for all $u_1,\ldots, u_d \in[0,1]$, by
\[
\mathbb{\widetilde C}^{\m}_n (u_1, \ldots,
u_d)= \sqrt{n} \bigl[H^{\m} _n \bigl\{
F_{n1}^{\m-1} (u_1), \ldots, F_{nd}^{\m-1}
(u_d) \bigr\} -
H^{\m} \bigl\{ F_{n1}^{\m-1} (u_1),
\ldots, F_{nd}^{\m-1} (u_d) \bigr\} \bigr]
\]
and
\[
\mathbb{\widetilde D}_n(u_1,\ldots, u_d) =
\sqrt{n} \bigl[ C^{\m} \bigl\{ F_1^{\m}\circ
F_{n1}^{\m-1}(u_1),\ldots, F_d^{\m}
\circ F_{nd}^{\m
-1}(u_d) \bigr\} -
C^{\m}(u_1,\ldots, u_d) \bigr].
\]
The next proposition shows that $ \mathbb{\widetilde C}^{\m}_n$ has the
same asymptotic behavior as $\mathbb{C}^{\m}_n$.
%
%
%pr4.2 #&#
\begin{prop}\label{prop42}
As $n \to\infty$, $\|\mathbb{C}_{n}^{\m}- \mathbb{\widetilde
C}_{n}^{\m}
\| \stackrel{\mathrm{ p}}{\to} 0$.
\end{prop}

Next, one needs to determine the limit of the second summand in \eqref
{eq41}. The following result first shows that $\mathbb{\widetilde
D}_n$ has the same asymptotic behavior as that of the auxiliary process
$\mathbb{D}_n$ defined, for all $u_1,\ldots, u_d \in[0,1]$, by
\[
\mathbb{D}_n(u_1,\ldots, u_d) = \sqrt{n}
\biggl[C^{\m} \biggl\{ u_1 - \frac
{\mathbb{C}_{n1}^{\m}(u_1)}{\sqrt{n}} , \ldots,
u_d - \frac{\mathbb
{C}_{nd}^{\m}(u_d)}{\sqrt{n}} \biggr\} - C^{\m}(u_1,
\ldots, u_d) \biggr],
\]
where $\mathbb{C}_{n1}^{\m}, \ldots,\mathbb{C}_{nd}^{\m}$ are the
margins of
$\mathbb{C}^{\m}_n$.

%
%pr4.3 #&#
\begin{prop}\label{prop43}
As $n \to\infty$, $\|\mathbb{D}_{n} - \mathbb{\widetilde D}_{n} \|
\stackrel{\mathrm{ p}}{\to} 0$.
\end{prop}

Finally, fix an arbitrary compact subset $K$ of $\mathcal{O}$ and
consider the mapping $\mathfrak{D}_K \dvt\break \mathcal{C}([0,1]^d) \to
\mathcal
{C}(K)$ defined, for all $g\in\mathcal{C}([0,1]^d)$ and $(u_1,\ldots,
u_d) \in K$, by
\[
\mathfrak{D}_K(g) (u_1,\ldots, u_d) = - \sum
_{j=1}^d {\dot C} ^{\m}
_j(u_1,\ldots, u_d) g(1,\ldots, 1,
u_j, 1,\ldots, 1).
\]
This mapping is clearly linear and continuous because for any $g$, $g^*
\in\mathcal{C}([0,1]^d)$,
\[
\bigl\| \mathfrak{D}_K(g) - \mathfrak{D}_K \bigl(g^* \bigr)\bigr\|
\le\sum_{j=1}^d {\dot C}^{\m}_j(u_1,
\ldots, u_d) \bigl\| g - g^* \bigr\| \le d \bigl\| g -g^*\bigr\|.
\]
For, when they exist, the partial derivatives of any copula take values
in $[0,1]$. The Continuous Mapping theorem then implies that, as $n\to
\infty$, $\mathfrak{D}_K(\mathbb{C}_n^{\m}) \rightsquigarrow
\mathfrak
{D}_K(\mathbb{C}^{\m})$ in $\mathcal{C}(K)$.
As shown next, the difference between $\mathbb{ D}_n$ and $\mathfrak
{D}_K(\mathbb{C}_n^{\m}) $ is asymptotically negligible.

%
%pr4.4 #&#
\begin{prop}\label{prop44}
As $n \to\infty$, $\| \mathbb{D}_{n} - \mathfrak{D}_K (\mathbb
{C}_n^{\m}
) \|_K \stackrel{\mathrm{ p}}{\to} 0$ for any compact $K \subset\mathcal{O}$.
\end{prop}

\com{To complete the proof of Theorem~\ref{thm31}, let $K$ be any
compact subset of $\mathcal{O}$. Combining Propositions \ref
{prop41}--\ref{prop44}, one finds that, as $n \to\infty$,
\[
\bigl\| \mathbb{\widehat C}_n^{\scriptsize\maltese}- \mathbb{C}_n^{\scriptsize\maltese}-
\mathfrak{D}_K \bigl(\mathbb{C}_n^{\m} \bigr)
\bigr\|_K \stackrel{\mathrm{ p}} {\to} 0.
\]
The Continuous Mapping theorem can then be invoked together with
Proposition \ref{prop41} to conclude that $\mathbb{\widehat
C}^{\scriptsize\maltese}_n \rightsquigarrow\mathbb{\widehat C}^{\scriptsize\maltese}=
\mathbb
{C}^{\scriptsize\maltese}+ \mathfrak{D}_K (\mathbb{C}^{\m})$. To establish Corollary
\ref{cor31}, first note that when $C^{\m}= \Pi$, $\dot C_j^{\m}$ is
continuous on $[0,1]^d$ for all $j\in\{1,\ldots,d\}$. One can then
define $\mathfrak{D}$ as $\mathfrak{D}_K$ with $K=[0,1]^d$ and use the
following result to conclude.}
\com{
%
%
%pr4.5 #&#
\begin{prop}\label{prop45}
When $C^{\m}= \Pi$, $\| \mathbb{D}_{n} - \mathfrak{D}(\mathbb
{C}_n^{\m})
\| \stackrel{\mathrm{ p}}{\to} 0$ as $n\to\infty$.
\end{prop}
}

%
%re4.1 #&#
\begin{rem}\label{rem41}
Although the process $\mathbb{\widehat C}_n^{\m}$ fails to converge on
$\mathcal{C}([0,1]^d)$ in general, the sequence $\| \mathbb{\widehat
C}_n^{\m}\|$ is tight. Indeed, the definition of $\mathbb{D}_n$ and the
Lipschitz property of $C^{\m}$ imply that $ \| \mathbb{D}_n \| \le\|
\mathbb{C}^{\m}_{n1} \| + \cdots+ \| \mathbb{C}^{\m}_{nd} \| \le d
\|
\mathbb{C}_n^{\m}\|$. From $\eqref{eq41}$ and the triangle inequality,
%
%
%e4.2 #&#
\begin{equation}
\label{eq42} \bigl\| \mathbb{\widehat C}_n^{\m}\bigr\| \le(d+1) \bigl\|
\mathbb{C}_n^{\m}\bigr\| +\bigl \| \mathbb{C}_n^{\m}-
\mathbb{\widetilde C}_n^{\m}\bigr\| + \| \mathbb{D}_n
- \mathbb{\widetilde D}_n \|.
\end{equation}
The result thus follows because the three summands form tight
sequences. Indeed, $\mathbb{C}_n^{\m}$ converges weakly in $\mathcal
{C}([0,1]^d)$ by Proposition~\ref{prop41} and the other two terms
converge in probability to $0$ by Propositions~\ref{prop42} and \ref
{prop43}, respectively. It is further of interest to observe that
because $\|\dot C^{\m}_j\|_{\mathcal{O}} \le1$ for all $j \in\{ 1,
\ldots, d\}$, one has $\|\mathbb{\widehat C}^{\m}\|_{\mathcal{O}}
\le
(d+1)\|\mathbb{C}^{\m}\|$.
\end{rem}

Finally, note that $\widehat C_n^{\m}$ is a uniformly consistent
estimator of $C^{\m}$. This follows immediately from \eqref{eq42}, the
Continuous Mapping theorem and Slutsky's lemma.

%
%co4.1 #&#
\begin{cor}\label{cor41}
As $n\to\infty$, $\| \widehat C_n^{\m}- C^{\m}\| \stackrel{\mathrm{ p}}{\to} 0$.
\end{cor}

%s5 #&#
\section{Applications} \label{sec5}

Theorem~\ref{thm31} characterizes the weak limit of the empirical
process $\mathbb{\widehat C}^{\m}_n$ in $\mathcal{C}(K)$ for any compact
subset $K$ of $\mathcal{O}$. \com{To illustrate the usefulness of this
result for inference, a few initial examples of application are
provided below. They pertain to classical statistics for monotone trend
and tests of independence, respectively.}

%s5.1 #&#
\subsection{Tests of monotone trend}

Kendall's tau and Spearman's rho are two classical measures of monotone
trend for two-way cross-classifications of ordinal or interval data. As
described, for example, in Agresti \cite{Agresti2007}, powerful tests of
independence can be based on these statistics. Both of them are
functions of (mid-) ranks that can be expressed as functionals of
$\widehat C_n^{\m}$ (Ne{\v{s}}lehov{\'a} \cite{Neslehova2007}).

Given a random sample $\mathcal{X} = \{ (X_{11}, X_{12}), \ldots,
(X_{n1}, X_{n2}) \}$ from a bivariate distribution function $H$, let
$R_{ij}$ denote the component-wise mid-rank of $X_{ij}$ for $i \in\{
1,\ldots,n\}$ and $j\in\{1,2\}$. Let also $a_n$ and $b_n$, respectively,
represent the number of strictly concordant and discordant pairs in the
sample. The non-normalized versions of Kendall's and Spearman's
coefficients then satisfy
\begin{eqnarray*}
\tau_n & =& \frac{a_n - b_n}{{n \choose2}} = \frac{n-1}{n} \biggl\{
- 1 + 4
\int_0^1 \int_0^1
{\widehat C}_n^{\m}(u,v) \,\d{\widehat C}_n^{\m}(u,v)
\biggr\},
\\
\rho_n & =& \frac{12}{n^3} \sum_{i=1}^n
\biggl(R_{i1} - \frac{n+1}{2} \biggr) \biggl( R_{i2} -
\frac{n+1}{2} \biggr) = 12 \int_{[0,1]^2} \bigl\{ \widehat
C^{\m}_n (u,v) - uv \bigr\} \,\d\Pi(u,v).
\end{eqnarray*}
It is immediate from Corollary \ref{cor41} that $\tau_n$ and $\rho_n$
are consistent estimators of
\[
\tau= - 1 + 4 \int_0^1 \int
_0^1 C^{\m}(u,v) \,\d
C^{\m}(u,v), \qquad\rho= 12 \int_{[0,1]^2} \bigl\{
C^{\m}(u,v) - uv \bigr\} \,\d\Pi(u,v).
\]
It is well known that $\tau_n$ is a $U$-statistic and hence
asymptotically Gaussian (Lee \cite{Lee1990}). Its limiting behavior can
also be deduced from Theorem \ref{thm31}. To see this, first call on
Hoeffding's identity (Nelsen \cite{Nelsen1999}, Corollary 5.1.2) to write
\[
\int_{[0,1]^2} C^{\m}(u,v) \,\d{\widehat
C}_n^{\m}(u,v) = \int_{[0,1]^2} {\widehat
C}_n^{\m}(u,v) \,\d C^{\m}(u,v).
\]
Given that $\widehat C_n^{\scriptsize\maltese}$ and $C^{\scriptsize\maltese}$ are absolutely
continuous with respect to the Lebesgue measure, the fact that the
complement of $\mathcal{O}$ in $[0,1]^d$ has Lebesgue measure $0$ then
implies that
\begin{eqnarray*}
&&\sqrt{n} \biggl\{ \int_{[0,1]^2} {\widehat C}_n^{\m}(u,v)
\,\d{\widehat C}_n^{\m}(u,v) - \int_{[0,1]^2}
C^{\m}(u,v) \,\d C^{\m}(u,v) \biggr\}
\\
&&\quad= \int_{\mathcal{O}} \mathbb{\widehat C}_n^{\m}(u,v)
\,\d{\widehat C}_n^{\m} (u,v) + \int_{\mathcal{O}}
\mathbb{\widehat C}_n^{\m}(u,v) \,\d C^{\m}(u,v)
.
\end{eqnarray*}

The following representation for the limit of $\sqrt{n} (\tau_n -
\tau)$ can be deduced from this relation. Details are provided in
Appendix \ref{appC}.

%
%pr5.1 #&#
\begin{prop}
\label{prop51}
In dimension $d = 2$, $\sqrt{n} (\tau_n - \tau)$ converges weakly,
as $n \to\infty$, to the centred Gaussian random variable
\[
\mathcal{T}_2 = 8 \int_{\mathcal{O}} \mathbb{\widehat
C}^{\m}(u,v) \,\d C^{\m}(u,v).
\]
\end{prop}

Similarly, the asymptotic normality of $\sqrt{n} (\rho_n -\rho)$ can
be deduced from the theory of $U$-statistics; see, for example, Quessy
\cite
{Quessy2009b}. The latter paper also considers several $d$-variate
extensions of $\rho_n$ which mimic the multivariate versions of these
coefficients for continuous data proposed by Schmid and
Schmidt \cite
{SchmidSchmidt2007}. Recently, we proposed alternative estimators of
$\rho$ in the multivariate case and showed that they lead to powerful
tests of independence and a graphical tool for visualizing dependence
in discrete data (Genest, Ne{\v{s}}lehov\'a and
R\'emillard \cite{GenestNeslehovaRemillard2013}). In
particular, we considered
\[
\rho_{nd} = \varrho_d \Biggl[ -\frac{1}{2^d} +
\frac{1}{n}\sum_{i=1}^n \Biggl\{
\prod_{j=1}^d \biggl( \frac{2n+1}{2n} -
\frac{R_{ij}}{n} \biggr) \Biggr\} \Biggr],
\]
where $\varrho_d = 2^d(d+1)/\{2^d-(d+1)\}$ and for each $i \in\{ 1,
\ldots, n \}$ and $j \in\{ 1, \ldots, d\}$, $R_{ij}$ denotes the
mid-rank of $X_{ij}$ among $X_{1j}, \ldots, X_{nj}$. The latter
reduces to $\rho_n$ in the bivariate case and can be rewritten as
\[
\rho_{nd} = \varrho_d \int_{[0,1]^d} \bigl
\{ \widehat C^{\m}_n (u_1,\ldots,
u_d) - \Pi(u_1,\ldots, u_d) \bigr\} \,\d\Pi
(u_1,\ldots, u_d).
\]
Furthermore, it is a consistent estimator of
\[
\rho_d = \varrho_d \int_{[0,1]^d} \bigl\{
C^{\m}(u_1, \ldots, u_d) - \Pi
(u_1,\ldots, u_d) \bigr\} \,\d\Pi(u_1,
\ldots, u_d).
\]
The asymptotic normality of $\sqrt{n} (\rho_{nd} - \rho_d)$,
established by Genest, Ne{\v{s}}lehov\'a and
R\'emillard~\cite{GenestNeslehovaRemillard2013}, can be shown
alternatively using Theorem \ref{thm31}. A detailed proof of the
following result is given in Appendix \ref{appC}.
%
%
%pr5.2 #&#
\begin{prop}
\label{prop52}
In arbitrary dimension $d \ge2$, $\sqrt{n} (\rho_{nd} - \rho_d)$
converges weakly, as $n \to\infty$, to the centred Gaussian random variable
\[
\mathcal{R}_d = \varrho_d \int_{\mathcal{O}}
\mathbb{\widehat C}^{\m} (u_1,\ldots, u_d) \,\d
\Pi(u_1,\ldots, u_d) .
\]
\end{prop}

%s5.2 #&#
\subsection{Tests of independence}

\com{When dealing with contingency tables that are sparse or whose
dimension varies with the sample size, Theorem \ref{thm31} can be
used to construct consistent and powerful tests of independence. This
is because random variables $X_1, \ldots, X_d$ are mutually
independent if and only if $C^{\scriptsize\maltese}= \Pi$.} To test the null
hypothesis $\mathcal{H}_0$ of mutual independence between $X_1, \ldots
, X_d$, one could consider, for example, the Cram\'er--von Mises statistic
\[
S_n= n\int_{[0,1]^d} \bigl\{\widehat
C_n^{\m}(u_1,\ldots, u_d) - \Pi
(u_1,\ldots, u_d) \bigr\}^2 \,\d
\Pi(u_1,\ldots, u_d).
\]
Note that when $X_1,\ldots, X_d$ are continuous, $S_n$ is equivalent to
the statistic suggested by Deheuvels \cite{Deheuvels-Raoult1979} and later
studied by Genest and
R{\'e}millard \cite{GenestRemillard2004}. The limiting distribution of
$S_n$ under $\mathcal{H}_0$ is easily deduced from Corollary \ref
{cor31} when the variables are integer-valued or increasing
transformations thereof. In fact, a straightforward adaptation of the
proof of Proposition \ref{prop52} yields the following result.

%
%pr5.3 #&#
\begin{prop}\label{prop53} Under $\mathcal{H}_0$ one has, as $n\to
\infty$, $S_n \rightsquigarrow S$, where
\[
S= \int_{[0,1]^d} \bigl\{ \widehat{\mathbb{C}}^{\m}(u_1,
\ldots, u_d) \bigr\}^2 \,\d\Pi(u_1,\ldots,
u_d).
\]
If $\mathcal{H}_0$ does not hold, then, as $n\to\infty$,
\[
\frac{S_n}{n} \stackrel{\mathrm{ p}} {\to} \int_{[0,1]^d} \bigl\{
C^{\m}(u_1,\ldots, u_d) - \Pi(u_1,
\ldots, u_d) \bigr\}^2 \,\d\Pi(u_1,\ldots,
u_d) >0.
\]
\end{prop}

In particular, Proposition \ref{prop53} implies that a test based on
$S_n$ is consistent against any alternative, that is, when $\mathcal
{H}_0$ fails then, as $n\to\infty$, $\Pr(S_n>\varepsilon)\to1$ for all
$\varepsilon>0$.

\com{Unfortunately, the limiting null distribution of $S_n$ depends on
the margins of $H$ which are generally unknown. To carry out the test,
one must thus resort to resampling techniques, such as the multiplier
bootstrap (van~der Vaart and
Wellner \cite{VaartWellner1996}). An illustration of how this can
be done is presented below in the case $d=2$.}

\com{
%
%
%al5.1 #&#
\begin{algo}
Given a random sample $\mathcal{X} = \{ (X_{11}, X_{12}), \ldots,
(X_{n1}, X_{n2}) \}$ from a bivariate distribution function $H$,
define, for $i\in\{1,\ldots,n\}$ and $j\in\{1,2\}$,
\[
V_{nj,i}(u) = \lambda_{F_{nj}}(u)\I\bigl\{X_{ij}\le
A_{nj}(k_j) \bigr\} + \bigl\{1- \lambda_{F_{nj}}(u)
\bigr\}\I\bigl\{X_{ij}\le A_{nj}(k_j-1) \bigr\},
\]
whenever $F_{nj}\{A_{nj}(k_j-1)\} < u \le F_{nj}\{A_{nj}(k_j)\}$ for
some $k_j \in\{0,\ldots, n_j\}$. The test based on $S_n$ can now be
carried out as follows.
\begin{enumerate}
\item[\emph{Step} 1:] For each $m \in\{1,\ldots, M\}$, generate an independent
random sample $\xi^{(m)}_1,\ldots, \xi^{(m)}_n$ of size $n$ from a
univariate distribution with mean zero and variance 1, and set $\bar
\xi
^{(m)} = (\xi^{(m)}_1+\cdots+ \xi^{(m)}_n)/n$.
\item[\emph{Step} 2:] For each $m \in\{1,\ldots, M\}$, define the process
$\mathfrak{C}_n^{(m)}$ at each $u,v \in[0,1]$ by
\[
\mathfrak{C}_n^{(m)}(u,v) = \frac{1}{\sqrt{n}}\sum
_{i=1}^n \bigl(\xi_i^{(m)} -
\bar\xi^{(m)} \bigr) \bigl\{V_{n1,i}(u) - u \bigr\} \bigl
\{V_{n2,i}(v) - v \bigr\}
\]
and compute
\[
S_n^{(m)} = \int_0^1
\int_0^1 \bigl\{\mathfrak{C}_n^{(m)}(u,v)
\bigr\}^2 \,\mathrm{d} v \,\mathrm{d} u.
\]
\item[\emph{Step} 3:] Estimate the $p$-value for the test by
\[
\frac{1}{M} \sum_{m=1}^M \I
\bigl(S_n^{(m)} > S_n \bigr).
\]
\end{enumerate}
\end{algo}}

An efficient implementation of this procedure is described in a
companion paper in preparation, in which the
validity of the multiplier bootstrap is established in this specific
context. Here, the finite-sample properties of this test are merely
illustrated through a small simulation study involving:
\begin{itemize}
\item five copulas: independence, Clayton (Cl) and Gaussian (Ga) with
$\tau\in\{ 0.1, 0.2\}$;
\item four margins: $\operatorname{Binomial}(3, 0.5)$, $\operatorname{Poisson}(1)$, $\operatorname{Poisson}(20)$,
$\operatorname{Geometric}(0.5)$, respectively, denoted by $F_1$, $F_2$, $F_3$ and $F_4$;
\item three statistics: $S_n$, the standard $\chi^2$, and a modified
version available in \textsf{R} in which the $p$-value is computed by a
Monte Carlo method;
\item sample size $n = 100$ and nominal level $\alpha= 5\%$;
\item$M=1000$ multiplier replicates and $N=1000$ repetitions of the simulation.
\end{itemize}
The results of the study are displayed in Table~\ref{tab51} below.
The test based on $S_n$ maintains its nominal level very well in every
scenario. In contrast, the standard $\chi^2$ statistic performs rather
poorly except when one of the margins is $F_1$. Resorting to the Monte
Carlo $\chi^2$ statistic improves the level, but the test is still
slightly liberal in some cases.

The power of the test based on $S_n$ is way better than that of its two
competitors in columns~1--7 and~9. In columns~8 and~10, $\chi^2$ is
slightly better when $\tau=0.1$. Note however that in these cases, the
level of the $\chi^2$ statistic is completely off. For a more thorough
simulation study, see Murphy~\cite{Murphy2013}.

%
%t1 #&#
\begin{table}
\caption{Percentage of rejection of the null hypothesis $\mathcal{H}_0$
of mutual independence for the three tests considered in the simulation
study under various conditions}\label{tab51}
\begin{tabular*}{\textwidth}{@{\extracolsep{\fill}}llld{2.1}d{2.1}d{2.1}d{2.1}d{2.1}d{2.1}d{2.1}d{2.1}d{2.1}d{2.1}@{}}
\hline
& &&\multicolumn{10}{l@{}}{Distribution of $X_1$}\\[-6pt]
&&& \multicolumn{10}{l@{}}{\hrulefill}\\
& &&\multicolumn{1}{l}{$F_1$} & \multicolumn{1}{l}{$F_1$} &
\multicolumn{1}{l}{$F_2$} & \multicolumn{1}{l}{$F_1$} &
\multicolumn{1}{l}{$F_2$} & \multicolumn{1}{l}{$F_3$} &
\multicolumn{1}{l}{$F_1$} & \multicolumn{1}{l}{$F_2$} &
\multicolumn{1}{l}{$F_3$} &\multicolumn{1}{l@{}}{$F_4$}\\[-6pt]
&&& \multicolumn{10}{l@{}}{\hrulefill}\\
&&& \multicolumn{10}{l@{}}{Distribution of $X_2$}\\[-6pt]
&&& \multicolumn{10}{l@{}}{\hrulefill}\\
 &&& \multicolumn{1}{l}{$F_1$} & \multicolumn{1}{l}{$F_2$} &
 \multicolumn{1}{l}{$F_2$} & \multicolumn{1}{l}{$F_3$} &
 \multicolumn{1}{l}{$F_3$} & \multicolumn{1}{l}{$F_3$} &
\multicolumn{1}{l}{$F_4$} & \multicolumn{1}{l}{$F_4$} &
\multicolumn{1}{l}{$F_4$} &\multicolumn{1}{l@{}}{$F_4$}\\
\hline
$\tau$ & $C$ & Test & & & & & & & & & & \\
\hline
0 & $\Pi$ & $S_n$ & 4.6 & 5.1 & 4.9 & 5.2 &5.0 & 4.8 & 4.5 & 4.8 & 4.9
& 5.0 \\
& & $\chi^2$ & 4.5 & 5.2 & 9.6 & 4.6 & 13.9 & 14.6 & 5.8 & 11.6 & 17.0
& 14.1\\
&& $\chi^2$-MC & 4.6 & 5.4 & 7.1 & 5.3 & 7.1 & 4.4 & 5.5 & 6.1 & 6.1 &
5.7\\[6pt]
$0.1$ & Cl & $S_n$ & 26.6 & 25.4 & 22.3 & 29.2 & 27.0 & 29.3 & 22.0
&19.1 & 22.5 & 18.2 \\
& & $\chi^2$ & 17.5 & 9.4 & 14.6 & 8.3& 10.2 & 23.7 & 6.2& 16.2& 11.7&
22.3 \\
& & $\chi^2$-MC & 17.6 & 9.8 &9.1 &9.9 &4.3 & 9.1& 6.5 &7.2& 3.6
&6.3\\[3pt]
& Ga & $S_n$ & 27.1 & 26.2 & 25.7 & 28.7 & 27.7 & 29.3 & 25.5 & 25.0 &
26.2 & 23.7 \\
& & $\chi^2$ &11.5 & 9.9 & 21.9 & 6.0 & 15.8 & 17.0 & 8.7 & 27.4 & 20.0
& 34.2 \\
& & $\chi^2$-MC & 11.9 & 10.6 & 15.7 & 6.5 & 7.7 & 7.2 & 8.4 & 14.5 &
7.8 & 13.1 \\[6pt]
$0.2$ & Cl & $S_n$ & 72.2 & 69.2 & 68.8 &76.6 &75.4 & 81.0 & 62.1 &
59.9 & 65.1 & 52.9 \\
& & $\chi^2$ & 56.3 & 34.1& 32.7 &27.6 &16.0& 42.6 &18.3 &32.0& 16.8&
36.5 \\
& & $\chi^2$-MC & 56.2 & 33.4 & 22.8 & 30.1 & 7.3 & 22.7 &18.4 &16.7 &
4.5 &13.0 \\[3pt]
& Ga & $S_n$ & 73.5 & 74.4 & 74.7 & 78.2 & 78.0 & 82.0 & 72.4 & 72.5 &
75.4 & 68.2 \\
& & $\chi^2$ & 41.9 & 33.8 & 52.1& 14.8& 29.8 &31.3 &25.3 &56.7& 33.7
&66.4 \\
& & $\chi^2$-MC & 43.2 & 33.4 & 39.5 &16.9 &15.2 &13.3& 26.7& 34.2&
14.1& 36.7\\
\hline
\end{tabular*}
\end{table}

%s6 #&#
\section{Conclusion}\label{sec6}

This paper considered the empirical multilinear copula process $\mathbb
{\widehat C}_n^{\m}$ based on count data. Its convergence was established
in $\mathcal{C}(K)$ for any compact $K\subset\mathcal{O}$, where
$\mathcal{O}$ is an open subset of $[0,1]^d$ avoiding the points at
which the first order partial derivatives of $C^{\m}$ do not exist. The
convergence of $\mathbb{\widehat C}_n^{\m}$ in $\mathcal{C}(K)$ is
sufficient to deduce the asymptotic behavior of simple functionals
thereof that are commonly used in statistical inference. This was
demonstrated in Section~\ref{sec5} using two standard measures of
association based on mid-ranks. While these specific results could have
been obtained using the theory of $U$-statistics, knowledge of the
limiting behavior of $\mathbb{\widehat C}_n^{\m}$ will be essential in
other situations. The new consistent test of independence studied in
Section~\ref{sec5} provides an example.

It is natural to ask whether the present findings can be extended to
the empirical multilinear copula process based on arbitrary
discontinuous data. Such an extension may well be possible, given that
the estimator $\widehat C_n^{\m}$ is defined in general. We are currently
investigating this issue. Once this task has been completed, the
process $\widehat{\mathbb C}_n^{\m}$ will provide a solid foundation for
inference in copula models with arbitrary margins.

%
%sA #&#
\begin{appendix}\label{app}
%%%%%%%%% Appendix %%%%%%%%%%%
%sA #&#
\section{Proofs from Section \texorpdfstring{\protect\ref{sec2}}{2}}

{\bf Proof of Proposition \ref{prop22}}. For all
$x_1,\ldots, x_d \in\mathbb{R}$, one has
\[
H^{\m}(x_1,\ldots, x_d) = \int
_{[0,1]^d} H(x_1+u_1,\ldots,
x_d+u_d) \,\d u_1 \cdots\,\d u_d.
\]
If $(x_1,\ldots, x_d) \in[k_1-1, k_1) \times\cdots\times
[k_d-1,k_d)$ for some $k_1,\ldots, k_d \in\mathbb{N}$, one can replace
each $x_j + u_j$ by $k_j - 1$ or by $k_j$, according as $0 < u_j < k_j
-x_j$ or $k_j -x_j \le u_j < 1$ because $H$ is supported on $\mathbb
{N}^d$. After straightforward simplification, it follows that
\[
H^{\m}(x_1,\ldots,x_d) = \sum
_{S\subseteq\{1,\ldots,d\}} H(k_S) \biggl\{ \prod
_{\ell\notin S} (k_\ell- x_\ell) \biggr\} \biggl\{
\prod_{\ell
\in S} (x_\ell- k_\ell+1)
\biggr\},
\]
where $k_S= (k_{S_1},\ldots, k_{S_d})$ and $k_{S_j} = k_j$ if $j \in S$
and $k_{S_j} = k_j -1$ otherwise.

If $F^{\m}$ is a generic margin of $H^{\m}$, then $F^{\m}$ is a linear
interpolation of $F$, that is, $F^{\m}(x) = 0$ for $x< -1$ while
%
%
%eA.1 #&#
\begin{equation}
\label{eqA1} F^{\m}(x) = F(k-1) + \Delta F(k) ( x- k + 1)
\end{equation}
when $x \in[k-1, k)$ for some $k \in\mathbb{N}$. Thus when $u \in
(F(k-1), F(k)]$, one has
%
%
%eA.2 #&#
\begin{equation}
\label{eqA2} F^{\m-1}(u) = k-1 + \frac{u-F(k-1)}{\Delta F(k)} .
\end{equation}
If $u=0$, one can set $F^{\m-1}(0)=-1$ for convenience, because the
support of $X$ is bounded below by $0$ by hypothesis. It is then
immediate that $H^{\m}\{F_1^{\m-1}(u_1),\ldots, F_d^{\m-1}(u_d)\}$
yields the formula for $C^{\m}$ given in Proposition \ref{prop21}.
\hfill$\Box$

%sA #&#
\section{Proofs from Section \texorpdfstring{\protect\ref{sec4}}{4}}\label{appB}

The following elementary result is used in the sequel.

%
%leA.1 #&#
\begin{lem}
\label{lemB1} If $G$ is a cumulative distribution function, then for
all $u\in(0,1)$ and $x\in\mathbb{R}$, one has $u \le G(x)
\Leftrightarrow G^{-1}(u) \le x \Leftrightarrow G \circ G^{-1} (u) \le G(x)$.
\end{lem}

\begin{pf*}{Proof of Proposition \ref{prop41}} First note that
for fixed
values of $k_1, \ldots, k_d \in\mathbb{N}$, one has
\begin{eqnarray*}
D_n \bigl\{ F_1 (k_1), \ldots,
F_d (k_d) \bigr\} &=& \frac{1}{n} \sum
_{i=1}^n \prod_{j=1}^d
\I\bigl\{ F_j (X_{ij}) \le F_j
(k_j) \bigr\}
\\
&=& \frac
{1}{n} \sum_{i=1}^n \prod
_{j=1}^d \I\bigl\{ F_j \circ
F_j^{-1} (V_{ij}) \le F_j
(k_j) \bigr\}.
\end{eqnarray*}
In view of Lemma~\ref{lemB1}, it follows that
\[
D_n \bigl\{ F_1 (k_1), \ldots,
F_d (k_d) \bigr\} = \frac{1}{n} \sum
_{i=1}^n \prod_{j=1}^d
\I\bigl\{ V_{ij} \le F_j (k_j) \bigr\} =
B_n \bigl\{ F_1 (k_1), \ldots,
F_d (k_d) \bigr\} .
\]
From the definition of $\mathfrak{M}_H$, one then has $\mathfrak{M}_H
(D_n) = \mathfrak{M}_H (B_n)$ and hence $C_n^{\m}= \mathfrak{M}_H
(B_n)$. The linearity of $\mathfrak{M}_H$ and the fact that $\mathfrak
{M}_H (C^{\m}) = C^{\m}$ further imply that $\mathbb{C}^{\m}_n =
\mathfrak
{M}_H (\mathbb{B}_n)$ \com{from which it also follows that
%
%
%eA.1 #&#
\begin{equation}
\label{eqB1} \bigl\|\mathbb{C}^{\m}_n\bigr\| =\bigl \|
\mathfrak{M}_H (\mathbb{B}_n)\bigr\| \le\|\mathbb
{B}_n\|
\end{equation}
because the operator $\mathfrak{M}_H$ is a contraction.} Given that
$\mathfrak{M}_H$ is a continuous mapping and that $\mathbb
{B}_n\rightsquigarrow\mathbb{B}_{C^{\m}}$ as $n\to\infty$, the
Continuous Mapping theorem yields the conclusion.
\end{pf*}

The following auxiliary results are needed for the proof of Proposition
\ref{prop42}.

%
%leA.2 #&#
\begin{lem} \label{lemB2} For all $u_1, \ldots, u_d \in[0,1]$,
\[
C^{\m}_n (u_1, \ldots, u_d) =
H^{\m}_n \bigl\{ F_1^{\m-1}
(u_1), \ldots, F_d^{\m-1} (u_d) \bigr
\}.
\]
\end{lem}

\begin{pf}
First note that the functions on both sides of the above identity are
continuous on $[0,1]^d$. This is the case for $C^{\m}_n$, as explained in
Section~\ref{sec4}. To see why this is true for the other one, fix
arbitrary $u_1, \ldots, u_d \in[0,1)$ and observe that
\begin{eqnarray*}
&&\bigl| H^{\m}_n \bigl\{ F_1^{\m-1}
(u_1+), \ldots, F_d^{\m-1} (u_d+)
\bigr\} - H^{\m}_n \bigl\{ F_1^{\m-1}
(u_1), \ldots, F_d^{\m-1} (u_d) \bigr
\}\bigr|\\
&&\quad \le
\sum_{j=1}^d \bigl| F_{nj}^{\m}
\circ F_j^{\m-1}(u_j+) - F_{nj}^{\m}
\circ F^{\m-1}_j(u_j)\bigr|.
\end{eqnarray*}
Now each of the summands on the right-hand side must vanish. For, even
if $u_j$ is a point of discontinuity of $F_j^{\m-1}$ for some $j \in
\{ 1, \ldots, d \}$, the fact that $F_j^{\m}$ is continuous implies that
$F_j^{\m}\circ F_j^{\m-1}(u_j) = F_j^{\m}\circ F_j^{\m-1}(u_j+) =
u_j$. Now for arbitrary $x, y\in\mathbb{R}$, one has
%
%
%eA.2 #&#
\begin{equation}
\label{eqB2} F_j^{\m}(x) = F_j^{\m}(y)
\quad\Rightarrow \quad F_{nj}^{\m}(x) = F_{nj}^{\m}(y),
\end{equation}
because $F_{nj}$ can only jump where $F_j$ does. Hence $F_{nj}^{\m
}\circ
F_j^{\m-1} (u_j) = F_{nj}^{\m}\circ F_j^{\m-1}(u_j+)$.

Therefore, it suffices to look at the case where $u_1,\ldots, u_d \in
(0,1)$. Suppose that for each $j \in\{1 , \ldots, d \}$, $u_j \in
(F_j (k_j-1), F_j (k_j)]$ for some $k_j \in\mathbb{N}$. It then
follows from \eqref{eqA2} that, for all $j \in\{ 1, \ldots, d
\}$,
\[
F^{\m-1}_j (u_j) = k_j -1 +
\frac{u_j-F_j(k_j-1)}{\Delta F_j(k_j)} ,
\]
and hence
\begin{eqnarray*}
k_j- F^{\m-1}_j(u_j) &=&
\frac{F_j(k_j)-u_j}{\Delta F_j(k_j)} ,\\
 F^{\m-1}_j(u_j)
-k_j + 1 &=& \frac{u_j-F_j(k_j-1)}{\Delta F_j(k_j)} .
\end{eqnarray*}
Consequently,
\[
H^{\m}_n \bigl\{ F_1^{\m-1}
(u_1), \ldots, F_d^{\m-1} (u_d) \bigr
\} = \sum_{S
\subset\{ 1, \ldots, d\}} \lambda_{H,S}
(u_1,\ldots, u_d)H_n (k_{S_1},
\ldots, k_{S_d}).
\]
Now in view of Lemma \ref{lemB1}, one has
\begin{eqnarray*}
H_n(k_1,\ldots, k_d)& =& \frac{1}{n} \sum
_{i=1}^n \prod
_{j=1}^d \I(X_{ij} \le k_j) =
\frac{1}{n} \sum_{i=1}^n \prod
_{j=1}^d \I\bigl\{F_j^{-1}(V_{ij})
\le k_j \bigr\}
\\
&=& \frac{1}{n} \sum_{i=1}^n \prod
_{j=1}^d \I\bigl\{V_{ij} \le
F_j(k_j) \bigr\} = B_n \bigl\{
F_1(k_1),\ldots, F_d(k_d) \bigr\}.
\end{eqnarray*}
Therefore,
\[
H^{\m}_n \bigl\{ F_1^{\m-1}
(u_1), \ldots, F_d^{\m-1} (u_d) \bigr
\} = \sum_{S
\subset\{ 1, \ldots, d\}} \lambda_{H, S}
(u_1,\ldots, u_d) B_n \bigl\{
F_1(k_{S_1}), \ldots, F_d(k_{S_d})
\bigr\},
\]
which is $\mathfrak{M}_H(B_n)(u_1, \ldots, u_d)$. From the proof of
Proposition~\ref{prop41}, $\mathfrak{M}_H(B_n)=C^{\m}_n$.
\end{pf}

%
%leA.3 #&#
\begin{lem}\label{lemB3}For arbitrary $n \in\mathbb{N}$, $G_n =
F_n^{\scriptsize\maltese}\circ F^{\scriptsize\maltese-1}$ is a continuous distribution
function on $[0,1]$ and $G_n^{-1} = F^{\m}\circ F_n^{\m-1} $.
\end{lem}

\begin{pf}
As $G_n$ is the convolution of two non-decreasing functions, it is
non-decreasing. Furthermore, $G_n(0) = 0$ and $G_n(1) = 1$ by
construction. Proceeding as in the proof of Lemma \ref{lemB2}, one
can show that $G_n$ is indeed continuous. Turning to $G_n^{-1}$, fix
$u\in[0,1]$ and observe that for any $x\in\mathbb{R}$ such that
$F^{\m}
_n(x) \ge u$, one has $F^{\m}\circ F_n^{\m-1} (u) \le F^{\m}(x)$
because $F^{\m}$ is non-decreasing. Now suppose that $y \in\mathbb
{R}$ is
such that for all $x \in\mathbb{R}$, $F^{\m}_n(x) \ge u \Rightarrow y
\le F^{\m}(x)$. By virtue of Lemma \ref{lemB1}, this is equivalent to
saying that for all $x \in\mathbb{R}$, $F^{\m}_n(x) \ge u
\Rightarrow
F^{\m-1}(y) \le x$. This implies that $F^{\m-1}(y) \le F_n^{\m-1}
(u)$. Applying Lemma \ref{lemB1} once again, one can see that $y \le
F^{\m}\circ F_n^{\m-1} (u)$. Consequently,
\[
F^{\m}\circ F_n^{\m-1} (u) = \inf\bigl
\{F^{\m}(x) \dvt F^{\m}_n(x) \ge u \bigr\}.
\]

Next, $F^{\m}\circ F^{\m-1} (u) = u$ by continuity of $F^{\m}$. Hence,
for all $x\in\mathbb{R}$, $F^{\m}\circ F^{\m-1} \circ F^{\m}(x) =
F^{\m}
(x)$. Invoking implication \eqref{eqB2}, one deduces that $F^{\m}_n
\circ F^{\m-1} \circ F^{\m}(x) = F_n^{\m}(x)$, which implies
\begin{eqnarray*}
\inf\bigl\{F^{\m}(x) \dvt F^{\m}_n(x) \ge u
\bigr\} & = & \inf\bigl\{F^{\m}(x) \dvt F^{\m}_n
\circ F^{\m-1} \circ F^{\m}(x) \ge u \bigr\}
\\
& = & \inf\bigl\{v \dvt F^{\m}_n\circ F^{\m-1}(v)
\ge u \bigr\} = \inf\bigl\{v \dvt G_n (v) \ge u \bigr\}.
\end{eqnarray*}
In other words, $F^{\m}\circ F_n^{\m-1}= G_n^{-1}$.
\end{pf}

\begin{pf*}{Proof of Proposition \ref{prop42}} First note that
in view of Lemma \ref{lemB2} and Proposition \ref{prop22}, one has,
for all $u_1, \ldots, u_d \in[0,1]$,
\[
\mathbb{C}^{\m}_n (u_1, \ldots,
u_d) = \sqrt{n} \bigl[H^{\m}_n \bigl\{
F_1^{\m-1} (u_1), \ldots, F_d^{\m-1}
(u_d) \bigr\}
- H^{\m} \bigl\{ F_1^{\m-1} (u_1),
\ldots, F_d^{\m-1} (u_d) \bigr\} \bigr].
\]
Next observe that, for all $u_1,\ldots, u_d\in[0,1]$,
%
%
%eA.3 #&#
\begin{equation}
\label{eqB3} \mathbb{\widetilde C}_{n}^{\m}(u_1,
\ldots, u_d) = \mathbb{C}^{\m}_{n} \bigl\{
F_1^{\m}\circ F_{n1}^{\m-1}
(u_1), \ldots, F_d^{\m}\circ
F_{nd}^{\m-1} (u_d) \bigr\}.
\end{equation}
Indeed, one can write
\begin{eqnarray*}
&&\mathbb{C}^{\m}_{n} \bigl\{ F_1^{\m}
\circ F_{n1}^{\m-1} (u_1), \ldots,
F_d^{\m}\circ F_{nd}^{\m-1}
(u_d) \bigr\}\\
&&\quad =
\sqrt{n} \bigl[ H^{\m}_n \bigl\{F_1^{\m-1}
\circ F_1^{\m}\circ F_{n1}^{\m
-1}
(u_1), \ldots, F_d^{\m-1}\circ
F_d^{\m}\circ F_{nd}^{\m-1}
(u_d) \bigr\}
\\
&&\hspace*{16pt}\qquad{}- H^{\m} \bigl\{F_1^{\m-1} \circ
F_1^{\m}\circ F_{n1}^{\m-1}
(u_1), \ldots, F_d^{\m-1}\circ
F_d^{\m}\circ F_{nd}^{\m-1}
(u_d) \bigr\} \bigr].
\end{eqnarray*}
Furthermore,
\begin{eqnarray*}
&&\bigl| H^{\m}_n \bigl\{F_1^{\m-1} \circ
F_1^{\m}\circ F_{n1}^{\m-1}
(u_1), \ldots, F_d^{\m-1}\circ
F_d^{\m}\circ F_{nd}^{\m-1}
(u_d) \bigr\}
\\
&&\qquad{}- H^{\m}_n \bigl\{F_{n1}^{\m-1}
(u_1), \ldots, F_{nd}^{\m-1} (u_d) \bigr
\} \bigr|
\\
&&\quad\le\sum_{j=1}^{d} \bigl| F_{nj}^{\m}
\circ F_j^{\m-1}\circ F_j^{\m} \circ
F_{nj}^{\m-1} (u_j) - F_{nj}^{\m}
\circ F_{nj}^{\m-1} (u_j) \bigr| .
\end{eqnarray*}
Now the right-hand side is zero by Lemma \ref{lemB3} and the fact
that for all $j\in\{1,\ldots, d\}$ and $u_j \in[0,1]$, $F_{nj}^{\m}
\circ F_{nj}^{\m-1} (u_j) = u_j$ because $F_{nj}^{\m}$ is a continuous
distribution function. As $F_1^{\m}, \ldots, F^{\m}_d$ are also
continuous distribution functions, one has
\begin{eqnarray*}
&&\bigl| H^{\m} \bigl\{F_1^{\m-1} \circ
F_1^{\m}\circ F_{n1}^{\m-1}
(u_1), \ldots, F_d^{\m-1}\circ
F_d^{\m}\circ F_{nd}^{\m-1}
(u_d) \bigr\}
\\
&&\qquad{}- H^{\m} \bigl\{F_{n1}^{\m-1} (u_1),
\ldots, F_{nd}^{\m-1} (u_d) \bigr\}\bigr |
\\
&&\quad\le\sum_{j=1}^{d} \bigl| F_{j}^{\m}
\circ F_j^{\m-1}\circ F_j^{\m} \circ
F_{nj}^{\m-1} (u_j) - F_{j}^{\m}
\circ F_{nj}^{\m-1} (u_j) \bigr| =0.
\end{eqnarray*}
Therefore, identity \eqref{eqB3} holds and one can write
\[
\bigl\|\mathbb{C}^{\m}_{n} - \mathbb{\widetilde
C}_{n}^{\m}\bigr\| = \bigl\|\mathbb{C}^{\m}_{n} -
\mathbb{C}^{\m}_{n} \bigl\{F_1^{\scriptsize\maltese}
\circ F_{n1}^{\scriptsize\maltese-1},\ldots, F_d^{\scriptsize\maltese}\circ
F_{nd}^{\scriptsize\maltese-1} \bigr\}\bigr\|.
\]
\com{Next, using \eqref{eqA1} and \eqref{eqA2} applied to $F$ and
$F_n$, respectively, a direct calculation yields}
\begin{eqnarray*}
&&\sqrt{n} \bigl\{u_j -F_j^{\m}\circ
F_{nj}^{\scriptsize\maltese-1} (u_j) \bigr\}\\
&&\quad =
\mathbb{B}_{nj} \bigl\{F_j(k_j-1) \bigr\}
\biggl\{\frac{F_{nj}(k_j)-u_j} {\Delta
F_{nj}(k_j)} \biggr\} + \mathbb{B}_{nj} \bigl
\{F_j(k_j) \bigr\} \biggl\{\frac{u_j - F_{nj}(k_j-1)} {\Delta
F_{nj}(k_j)} \biggr\},
\end{eqnarray*}
whenever $u_j\in(F_{nj}(k_j-1), F_{nj}(k_j)]$ for some
$k_j\in\mathbb{N}$. It follows that
%
%
%eA.4 #&#
\begin{equation}
\label{eqB4} \sup_{u_j \in[0,1]} \bigl| F_j^{\m}\circ
F_{nj}^{\m-1} (u_j) - u_j \bigr| \le
\frac{1}{\sqrt{n}} \| \mathbb{B}_{nj} \| \le\frac{1}{\sqrt{n}} \|
\mathbb{B}_{n} \|.
\end{equation}
As $n \to\infty$, $\| \mathbb{B}_{n} \| \rightsquigarrow\| \mathbb
{B}_{C^{\m}} \|$ and hence $\| \mathbb{B}_{n} \| /\sqrt{n} \stackrel
{\mathrm{ p}}{\to} 0$. Now for arbitrary $\varepsilon>0$, one has
\[
P^* \bigl( \bigl\|\mathbb{C}^{\m}_{n} -\mathbb{\widetilde
C}^{\m}_{n}\bigr\| > \varepsilon\bigr) = P^* \bigl\{ \bigl\|
\mathbb{C}^{\m}_{n} -\mathbb{C}^{\m}_{n}
\bigl(F_1^{\scriptsize\maltese} \circ F_{n1}^{\m-1},
\ldots, F_d^{\m}\circ F_{nd}^{\m-1} \bigr)
\bigr\| > \varepsilon\bigr\},
\]
where $P^*$ denotes outer probability. Given $\delta> 0$, the
right-hand side is the same as
\begin{eqnarray*}
&&P^* \biggl\{ \bigl\|\mathbb{C}^{\m}_{n} -\mathbb{C}^{\m}_{n}
\bigl(F_1^{\scriptsize\maltese} \circ F_{n1}^{\m-1},
\ldots, F_d^{\m}\circ F_{nd}^{\m-1} \bigr)
\bigr\| > \varepsilon, \frac{\| \mathbb{B}_n\| } {\sqrt{n} } < \delta
\biggr\}
\\
&&\quad{}+ P^* \biggl\{ \bigl\|\mathbb{C}^{\m}_{n} -\mathbb{C}^{\m}_{n}
\bigl(F_1^{\m} \circ F_{n1}^{\m-1},
\ldots, F_d^{\m}\circ F_{nd}^{\m-1} \bigr)
\bigr\| > \varepsilon, \frac{\| \mathbb{B}_n\| } {\sqrt{n} } \ge\delta
\biggr\},
\end{eqnarray*}
and in view of \eqref{eqB4}, the latter is bounded above by
\[
P^* \bigl\{ \omega_n \bigl(\mathbb{C}_n^{\m},
\delta\bigr) > \varepsilon\bigr\} + P^* \biggl( \frac{\| \mathbb
{B}_n\| } {\sqrt{n} } \ge\delta
\biggr),
\]
where
\[
\omega_n \bigl(\mathbb{C}_n^{\m},\delta\bigr)
= \mathop{\sup_{u_j, v_j \in
[0,1]\dvt|u_j - v_j| < \delta, }}_{ j \in\{1,\ldots,
d\}}\bigl|\mathbb{C}^{\m}_{n}(u_1,
\ldots, u_d) -\mathbb{C}^{\m}_{n}
(v_1,\ldots, v_d)\bigr|.
\]
Therefore,
\[
\limsup_{n\to\infty} P^* \bigl( \bigl\|\mathbb{C}^{\m}_{n}
-\mathbb{\widetilde C}^{\m}_{n}\bigr\| > \varepsilon\bigr) \le
\limsup_{n\to\infty} P^* \bigl\{\omega_n \bigl(
\mathbb{C}_n^{\m},\delta\bigr) > \varepsilon\bigr\}.
\]

Finally, recall that $\mathbb{C}^{\m}_{n}$ converges weakly in
$\mathcal{C}([0,1]^d)$ to a measurable random element $\mathfrak{M}_H
(\mathbb{B}_{C^{\m}})$. Because $\mathcal{C}([0,1]^d)$ is complete and
separable, Theorem 11.5.4. in Dudley \cite{Dudley2002} implies that
$\mathfrak{M}_H (\mathbb{B}_{C^{\m}})$ is tight. It then follows from
Lemma 1.3.8. and Theorem 1.5.7. in van~der Vaart and
Wellner \cite{VaartWellner1996} that
the sequence $\mathbb{C}^{\m}_{n}$ is asymptotically tight and hence
asymptotically uniformly equicontinuous in probability, viz.
\[
\lim_{\delta\downarrow0} \limsup_{n\to\infty} P^* \bigl\{
\omega_n \bigl(\mathbb{C}_n^{\m},\delta\bigr) >
\varepsilon\bigr\} = 0.
\]
This means that as $n \to\infty$, $P^* (\|\mathbb{C}^{\m}_n -
\mathbb{\widetilde C}_n^{\m}\| > \varepsilon) \to0$ for all
$\epsilon> 0$.
\end{pf*}

\begin{pf*}{Proof of Proposition \ref{prop43}} For fixed $j
\in\{1,\ldots, d\}$ and $u_j \in[0,1]$, first write
$F^{\m}_j \circ F^{\m-1}_{nj} (u_j)$ in the form $u_j - \sqrt{n} \{
u_j - F^{\m}_j \circ F^{\m-1}_{nj} (u_j)\} / \sqrt{n}$. Then
\begin{eqnarray*}
\bigl\| \mathbb{D}_n - \mathbb{\widetilde D}_n \bigr\| &=& \sqrt{n}
\sup_{u_1,
\ldots, u_d \in[0,1]} \biggl| C^{\m} \biggl\{ u_1 -
\frac{\mathbb
{C}_{n1}^{\m}(u_1)}{\sqrt{n}} , \ldots, u_d - \frac{\mathbb
{C}_{nd}^{\m}(u_d)}{\sqrt{n}} \biggr\}
\\
&&\hspace*{45pt}\qquad{}- C^{\m}{ \biggl[} u_1 - \frac{\sqrt{n} \{ u_1 - F^{\m}_1 \circ
F^{\m-1}_{n1} (u_1)\}}{\sqrt{n}} , \ldots,
u_d \\
&&\hspace*{45pt}\qquad{}- \frac{\sqrt{n}
\{ u_d - F^{\m}_d \circ F^{\m-1}_{nd} (u_d)\}}{\sqrt{n}} { \biggr]}
\biggr| .
\end{eqnarray*}
The Lipschitz property of copulas further implies that
\[
\| \mathbb{D}_n - \mathbb{\widetilde D}_n \| \le\sum
_{j=1}^d \sup_{u_j\in[0,1]} \bigl
\llvert\mathbb{C}_{nj}^{\m}(u_j) - \sqrt{n}
\bigl\{ u_j - F^{\m}_j \circ
F^{\m-1}_{nj} (u_j) \bigr\} \bigr\rrvert.
\]
\com{For any $j\in\{1,\ldots, d\}$, one can now call upon Proposition
\ref{prop42} with $d=1$ and $H=F_j$
to conclude that, as $n\to\infty$, $\sup_{u_j\in[0,1]} \llvert
\mathbb{C}_{nj}^{\m}(u_j) - \sqrt{n} \{ u_j - F^{\m}_j \circ
F^{\m-1}_{nj} (u_j)\}\rrvert \stackrel{\mathrm{ p}}{\to} 0$.}
\end{pf*}

The proof of Proposition \ref{prop44} relies on the following lemma.

%
%leA.4 #&#
\begin{lem}\label{lemB4}
Let $u_1, \ldots, u_d \in[0,1]$ and $v_1,\ldots, v_d \in[0,1]$ be
such that for each $j \in\{ 1, \ldots, d\}$, $u_j, v_j \in(F_j
(k_j-1) , F_j (k_j))$ for some $k_j \in\mathbb{N}$. Then
\[
C^{\m}(v_1, \ldots, v_d) -
C^{\m}(u_1, \ldots, u_d) = \sum
_{m=1}^d (v_m - u_m) {\dot
C}^{\m}_m (w_{m1},\ldots, w_{md}),
\]
where $w_{mj}$ equals $u_j$ or $v_j$ according as $j < m$ or $j \ge
m$, respectively.
\end{lem}

\begin{pf} First, write $C^{\m}(v_1, \ldots, v_d) - C^{\m}(u_1,
\ldots, u_d)$ in the alternative form
\[
\sum_{m=1}^d \bigl\{ C^{\m}
(w_{m1},\ldots, w_{md}) - C^{\m} (w_{
(m+1)1},
\ldots, w_{(m+1)d}) \bigr\}.
\]
It must then be shown that for all $m \in\{1,\ldots, d\}$, one has
%
%
%eA.5 #&#
%eA.6 #&#
\begin{eqnarray}
\label{eqB5} &&{C}^{\m} (w_{m1},\ldots, w_{md}) -
{C}^{\m} (w_{(m+1)1},\ldots, w_{(m+1)d})
\nonumber
\\[-8pt]
\\[-8pt]
\nonumber
&&\quad= (v_m - u_m) {\dot C}^{\m}_m
(w_{m1},\ldots, w_{md}).
\end{eqnarray}

To this end, observe that on the left-hand side of \eqref{eqB5},
$C^{\m}
$ is evaluated at two vectors whose components are identical, except in
position $m$. Let $w_1, \ldots, w_{m-1}, w_{m+1}, \ldots, w_d$ be the
matching components and note that $w_{mm} = v_m$ while $w_{(m+1)m}
=u_m$. Given $S \subset\{ 1, \ldots, d\}$, let $s_m$ be the size of
$S \cap\{ m\}$. From the definition of $\lambda_{H,S}$, one has
\begin{eqnarray*}
\lambda_{H,S} (w_{m1},\ldots, w_{md}) &=&
\lambda_{H,S} (w_{1},\ldots, w_{m-1},
v_m, w_{m+1}, \ldots, w_{d})
\\
&= &\biggl\{ \frac{v_m - F_m(k_m-1)}{\Delta F_m (k_m)} \biggr\}^{s_m}
\times\biggl\{
\frac{F_m(k_m) - v_m}{\Delta F_m(k_m)} \biggr\}^{1-s_m}
\\
&&{}\times\biggl\{\mathop{\prod_{\ell\notin S}}_{ \ell\neq m}
\frac
{F_\ell(k_\ell) - w_\ell} {
\Delta F_\ell(k_\ell)} \biggr\} \times\biggl\{\mathop{\prod
_{\ell
\in S }}_{ \ell\neq m} \frac{w_\ell- F_\ell(k_\ell-1)} {\Delta
F_\ell(k_\ell)} \biggr\}
\end{eqnarray*}
and
\begin{eqnarray*}
&&\lambda_{H,S}(w_{(m+1)1},\ldots, w_{(m+1)d}) =
\lambda_{H,S} (w_{1},\ldots, w_{m-1},
u_m, w_{m+1}, \ldots, w_{d})
\\
&&\quad= \biggl\{ \frac{u_m -
F_m(k_m-1)}{\Delta F_m (k_m)} \biggr\}^{s_m} \times\biggl\{
\frac
{F_m(k_m) - u_m}{\Delta F_m (k_m)} \biggr\}^{1-s_m}
\\
&&\qquad{}\times\biggl\{\mathop{\prod_{\ell\notin S }}_{ \ell\neq m}
\frac{F_\ell(k_\ell) - w_\ell} {\Delta F_\ell(k_\ell)} \biggr\}
\times\biggl\{\mathop{\prod
_{\ell\in S}}_{ \ell\neq m} \frac{w_\ell- F_\ell(k_\ell-1)}
{\Delta F_\ell(k_\ell)} \biggr\}.
\end{eqnarray*}
Consequently, their difference is equal to
\[
(v_m - u_m) \frac{(-1)^{1-s_m}}{\Delta F_m (k_m)} \times\biggl\{
\mathop{
\prod_{\ell\notin S }}_{ \ell\neq m} \frac{F_\ell(k_\ell) -
w_\ell} {\Delta F_\ell(k_\ell)}
\biggr\} \times\biggl\{\mathop{\prod_{\ell\in S }}_{ \ell\neq m}
\frac{w_\ell- F_\ell(k_\ell-1)} {\Delta F_\ell(k_\ell)} \biggr\}.
\]
It then follows from the definition of $C^{\m}$ that
\begin{eqnarray*}
&&C^{\m}(w_{m1},\ldots, w_{md}) -
C^{\m}(w_{(m+1)1},\ldots, w_{(m+1)d})
\\
&&\quad= \sum_{S \subset\{ 1, \ldots, d\}} H(k_S) \bigl\{
\lambda_{H,S} (w_{m1},\ldots, w_{md}) -
\lambda_{H,S} (w_{(m+1)1},\ldots, w_{(m+1)d}) \bigr\}
\\
&&\quad= (v_m - u_m) \sum_{S \subset\{ 1, \ldots, d\}}
H(k_S) \frac
{(-1)^{1-s_m}}{\Delta F_m (k_m)}\biggl\{\mathop{\prod_{\ell\notin
S}}_{ \ell\neq m}
\frac{F_\ell(k_\ell) - w_\ell} {\Delta
F_\ell(k_\ell)} \biggr\} \times\biggl\{\mathop{\prod
_{\ell\in S }}_{
\ell\neq m} \frac{w_\ell- F_\ell(k_\ell-1)} {
\Delta F_\ell(k_\ell)} \biggr\}
\\
&&\quad= (v_m - u_m) {\dot C}^{\m}_m
(w_{m1},\ldots, w_{md}).
\end{eqnarray*}
This completes the argument.
\end{pf}

\begin{pf*}{Proof of Proposition \ref{prop44}} Recall from the
definition of $\mathcal{O}$ that because $K$ is compact, it can be
covered by finitely many open cubes of the form
\[
\mathcal{O}_\ell= \bigl(F_1(k_{1\ell}-1),
F_1(k_{\ell1}) \bigr) \times\cdots\times\bigl(F_d(k_{d\ell}
-1), F_d(k_{d\ell}) \bigr),
\]
where $k_{1\ell},\ldots, k_{d\ell} \in\mathbb{N}$ for $\ell\in
\{1,\ldots, L\}$. Given that the sets $\mathcal{O}_1 , \ldots,
\mathcal{O}_L$ are mutually disjoint, $K_\ell= K \cap
\mathcal{O}_\ell$ is compact for each $\ell\in\{1,\ldots, L\}$.
Therefore, $K = K_1 \cup\cdots\cup K_L$ is a union of finitely many
disjoint compact sets. For arbitrary $\delta> 0$, let
\[
K_{\ell, \delta} = \bigcup_{(x_1, \ldots, x_d) \in K_\ell} \bigl\{
(u_1, \ldots, u_d) \in\mathbb{R}^d \dvt|
u_1 - x_1| + \cdots+ |u_d - x_d| <
\delta\bigr\} .
\]
Because $K_1,\ldots, K_L$ are compact, there exists $\delta_0 > 0$
such that $K_{\ell,\delta_0} \subset\mathcal{O}_\ell$ for all
$\ell
\in\{1,\ldots, L\}$. Now fix $\delta^* < \delta_0$ and let $K^*$
denote the closure of $K_{\delta^*}=K_{1,\delta^*} \cup\cdots\cup
K_{d, \delta^*}$, which is compact. Then for all $\delta\in(0,
\delta^*)$, one has $K \subset K_{\delta} \subset K^* \subset
\mathcal{O}$. For fixed $\delta\in(0, \delta^*)$ and $\varepsilon>
0$, write
\begin{eqnarray*}
&&P^* \bigl\{\bigl\|\mathbb{D}_n - \mathfrak{D}_K \bigl(
\mathbb{C}_n^{\m} \bigr) \bigr\|_K > \varepsilon
\bigr\} \\
&&\quad= P^* \biggl\{ \bigl\|\mathbb{D}_n - \mathfrak{D}_K
\bigl( \mathbb{C}_n^{\m} \bigr)\bigr \|_K >
\varepsilon, \frac{\|\mathbb{B}_n\|}{\sqrt{n}} < \frac
{\delta}{d} \biggr\}
\\
&&\qquad{}+ P^* \biggl\{ \bigl\|\mathbb{D}_n - \mathfrak{D}_K \bigl(
\mathbb{C}_n^{\m} \bigr) \bigr\|_K > \varepsilon,
\frac{\|\mathbb{B}_n\|}{\sqrt{n}} \ge\frac{\delta}{d} \biggr\}.
\end{eqnarray*}
When the event $ \{ {\|\mathbb{B}_n\|}/{\sqrt{n}} <
{\delta}/{d} \}$ holds and
$(u_1,\ldots, u_d) \in K_\ell$ for some $\ell\in\{1,\ldots, L\}$,
\[
(v_1,\ldots, v_d) = \biggl(u_1 -
\frac{\mathbb{C}_{n1}^{\m}(u_1)}{\sqrt{n}}, \ldots, u_d - \frac
{\mathbb{C}_{nd}^{\m}(u_d)}{\sqrt{n}} \biggr) \in
K_{\ell,
\delta}
\]
because, for all $j \in\{1,\ldots, d\}$, $\|
\mathbb{C}_{nj}^{\m}\| \le\| \mathbb{C}_{n}^{\m}\| \le\|
\mathbb{B}_n\|$ by \eqref{eqB1}. From Lemma \ref{lemB4},
\begin{eqnarray*}
&&\bigl\llvert\mathbb{D}_n(u_1,\ldots,
u_d) - \mathfrak{D}_K \bigl(\mathbb{C}_n^{\m}
\bigr) (u_1,\ldots, u_d) \bigr\rrvert
\\
&&\quad= \Biggl| \sum_{j=1}^d \mathbb{C}_{nj}^{\m}(u_j)
\bigl\{ {\dot C}_j^{\m} (u_1,\ldots,
u_d) - {\dot C}_j^{\m}(u_1,\ldots,
u_{j-1}, v_j, \ldots, v_d) \bigr\} \Biggr|
\\
&&\quad \le \| \mathbb{B}_n\| \sum_{j=1}^d
\bigl|  {\dot C}_j^{\m} (u_1,\ldots,
u_d) - {\dot C}_j^{\m} (u_1,\ldots,
u_{j-1}, v_j, \ldots, v_d )  \bigr|.
\end{eqnarray*}
Consequently, $\|\mathbb{D}_n - \mathfrak{D}_K(\mathbb{C}_n^{\m})
\|_K$ is bounded above by
\begin{eqnarray*}
&&\| \mathbb{B}_n\| \sum_{j=1}^d
\sup_{(u_1,\ldots, u_d) \in
K} \bigl\llvert {\dot C}_j^{\m}(u_1,
\ldots, u_d) - {\dot C}_j^{\m}
(u_1, \ldots, u_{j-1}, v_j, \ldots,
v_d )  \bigr\rrvert
\\
&&\quad \le \| \mathbb{B}_n\| \sum_{j=1}^d
\mathop{ \mathop{\sup_{(u_1,\ldots, u_d)
\in K, }}_{
(w_1,\ldots, w_d) \in K_\delta, }}_{
\sum_{m=1}^d | u_m - w_m| < \delta}
\bigl\llvert{\dot C}_j^{\m}(u_1,\ldots,
u_d) - {\dot C}_j^{\m} (w_1,\ldots,
w_d )  \bigr\rrvert
\\
&&\quad \le \| \mathbb{B}_n\| \omega(\delta),
\end{eqnarray*}
where\vspace*{1pt}
\[
\omega(\delta) = \sum_{j=1}^d \mathop{
\sup_{(u_1,\ldots,
u_d),(w_1,\ldots, w_d) \in K^*,}}_{
\sum_{m=1}^d | u_m - w_m| < \delta} \bigl\llvert {\dot
C}_j^{\m} (u_1,\ldots, u_d) - {\dot
C}_j^{\m}(w_1,\ldots, w_d) \bigr\rrvert.
\]
This observation implies that
\begin{eqnarray*}
&&\limsup_{n\to\infty} P^* \biggl\{ \bigl\|\mathbb{D}_n
- \mathfrak{D}_K \bigl(\mathbb{C}_n^{\m}
\bigr) \bigr\|_K > \varepsilon, \frac{\|\mathbb{B}_n\|}{\sqrt{n}} <
\frac{\delta}{d}
\biggr\}
\\
&&\quad\le \limsup_{n\to\infty} P^* \biggl\{ \| \mathbb{B}_n\|
\omega(\delta) > \varepsilon, \frac{\|\mathbb{B}_n\|}{\sqrt{n}} <
\frac{\delta
}{d} \biggr\}
\\
&&\quad\le \limsup_{n\to\infty} P^* \bigl\{ \| \mathbb{B}_n\|
\omega(\delta) > \varepsilon\bigr\} = P^* \bigl\{ \| \mathbb
{B}_{C^{\m}}\|
\omega(\delta) > \varepsilon\bigr\},
\end{eqnarray*}
where the equality is justified by the fact that $\| \mathbb{B}_n\|
\rightsquigarrow\| \mathbb{B}_{C^{\m}}\|$ as $n\to\infty$. Now
$\omega
(\delta) \to0$ as $\delta\to0$ because ${\dot C}_j^{\m}$
is absolutely continuous on $K^*$ for all $j \in\{ 1,\ldots, d\}$.
Therefore, $P^* \{ \| \mathbb{B}_{C^{\m}}\| \omega(\delta) >
\varepsilon\} \to0$ as $\delta\to0$. Finally, observe that
\[
\limsup_{n\to\infty} P^* \biggl(\bigl\|\mathbb{D}_n -
\mathfrak{D}_K \bigl(\mathbb{C}_n^{\m} \bigr)
\bigr\|_K > \varepsilon, \frac{\|\mathbb{B}_n\|}{\sqrt{n}} \ge\frac
{\delta}{d} \biggr)
\le\limsup_{n\to\infty} P^* \biggl(\frac{\|\mathbb{B}_n\|}{\sqrt
{n}} \ge
\frac{\delta}{d} \biggr) = 0
\]
because $\|\mathbb{B}_n\|/\sqrt{n} \stackrel{\mathrm{ p}}{\to} 0$ as $n
\to\infty$. As $\varepsilon> 0$ is arbitrary, one can conclude.
\end{pf*}

The proof of Proposition \ref{prop45} relies on the following lemma.
\com{
%
%
%leA.5 #&#
\begin{lem}\label{lemB5} Let $G$ be the distribution function of a
uniform random variable on $(0,1)$. Then for every $j\in\{1,\ldots, d\}
$ and as $n\to\infty$,
\[
Y_{nj} = \sqrt{n} \sup_{0 \le u \le1} \biggl| G \biggl\{u-
\frac{\mathbb
{C}_{nj}^{\m}(u)}{\sqrt{n}} \biggr\} - u +\frac{\mathbb
{C}_{nj}^{\m}(u)}{\sqrt{n}} \biggr| \stackrel{\mathrm{p}} {\to} 0.
\]
\end{lem}
}

\begin{pf}
Fix $j \in\{1,\ldots, d\}$ and write
\begin{eqnarray*}
Y_{nj} &=& \sqrt{n} \sup_{0 \le u \le1} \biggl[ \biggl\{
\frac{\mathbb
{C}_{nj}^{\m}(u)}{\sqrt{n}} -u \biggr\}\I\biggl\{u < \frac
{\mathbb{C}_{nj}^{\m}
(u)}{\sqrt{n}} \biggr\}
\\
&&\hspace*{45pt}{}+ \biggl\{-\frac{\mathbb{C}_{nj}^{\m}(u)}{\sqrt{n}} -(1-u) \biggr
\}\I\biggl\{1-u < -
\frac{\mathbb{C}_{nj}^{\m}(u)}{\sqrt{n}} \biggr\} \biggr].
\end{eqnarray*}
Observe that if $\|\mathbb{C}_{nj}^{\m}\| \le M$ for some constant $M>0$,
then as $n \to\infty$,
\[
Y_{nj} \le\sup_{0 \le u \le M/\sqrt{n}}\bigl|\mathbb{C}_{nj}^{\m}(u)\bigr|
+ \sup_{1-M/\sqrt{n} \le u \le1}\bigl|\mathbb{C}_{nj}^{\m}(u)\bigr|
\stackrel{\mathrm{p}} {\to} 0
\]
because $\mathbb{C}_{nj}^{\m}\rightsquigarrow\mathbb{C}^{\m}_j$ in
$\mathcal{C}([0,1])$ and $\mathbb{C}^{\m}_j(0) = \mathbb{C}^{\m
}_j(1) =0$.
Now fix $\varepsilon> 0$ and invoke the tightness of $\mathbb
{C}_{nj}^{\m}$ to find $M>0$ such that $\Pr(\|\mathbb{C}_{nj}^{\m
}\| > M) <
\varepsilon/2$ for all $n \in\mathbb{N}$. Thus,
\begin{eqnarray*}
\Pr(Y_{nj} > \varepsilon) \le\Pr\bigl(\bigl\|\mathbb{C}_{nj}^{\m}
\bigr\| > M \bigr) + \Pr\Bigl( \sup_{0 \le u \le M/\sqrt{n}}\bigl|\mathbb
{C}_{nj}^{\m}(u)\bigr|
+ \sup_{1-M/\sqrt{n} \le u \le1}\bigl|\mathbb{C}_{nj}^{\m}(u)\bigr| >
\varepsilon\Bigr).
\end{eqnarray*}
If $n$ is large enough, the right-hand side of the above inequality is
at most $\varepsilon$.
\end{pf}

\begin{pf*}{Proof of Proposition \ref{prop45}} For all
$u_1,\ldots, u_d \in[0,1]$, let
\[
\mathbb{D}_n^*(u_1,\ldots,u_d) = \sqrt{n}
\Biggl[ \prod_{j=1}^d \biggl\{
u_j - \frac{\mathbb{C}_{nj}^{\m}(u)}{\sqrt{n}} \biggr\} - \prod
_{j=1}^d u_j \Biggr].
\]
Then
\[
\bigl\| \mathbb{D}_n - \mathbb{D}_n^* \bigr\| \le\sum
_{j=1}^d \sqrt{n} \sup_{0
\le u \le1} \biggl| G
\biggl\{u-\frac{\mathbb{C}_{nj}^{\m}(u)}{\sqrt{n}} \biggr\} - u
+\frac{\mathbb{C}_{nj}^{\m}(u)}{\sqrt{n}} \biggr|
\]
because $|\prod_{j=1}^d a_j - \prod_{j=1}^d b_j| \le\sum_{j=1}^d |a_j
- b_j|$ for all $a_1, \ldots, a_d$, $b_1, \ldots, b_d \in(0,1)$.
Lem\-ma~\ref{lemB5} thus implies that $\| \mathbb{D}_n - \mathbb{D}_n^*
\| \stackrel{\mathrm{p}}{\to} 0$ as $n\to\infty$. Now by the multinomial formula,
\[
\bigl\| \mathbb{D}_n^*-\mathfrak{D} \bigl(\mathbb{C}_n^{\scriptsize\maltese}
\bigr)\bigr\| \le\sqrt{n} \sum_{S\subset\{1,\ldots, d\}, |S| \ge2}
\prod
_{j \in S} \frac{\|
\mathbb{C}_{nj}^{\scriptsize\maltese}\|}{\sqrt{n}} \stackrel{\mathrm{p}} {\to} 0.
\]
\upqed\end{pf*}
%sA #&#
\section{Proofs from Section~\texorpdfstring{\protect\ref{sec5}}{5}}\label{appC}

The proofs of Propositions \ref{prop51} and \ref{prop52} have much
in common. \com{They both rely on the following straightforward
consequence of Proposition 6.3.9 in Brockwell and
Davis \cite{BrockwellDavis1991}.}

%
%leA.1 #&#
\begin{lem}\label{lemC1} Let $Z_n$ be a sequence of random variables.
Suppose that for all $\delta$, $\epsilon> 0$, there exists a sequence
$Y_{n,\delta, \epsilon}$ of random variables such that for all $n \in
\mathbb{N}$, $\Pr(|Z_n - Y_{n, \delta, \epsilon}| > \delta) <
\epsilon
$ and $Y_{n,\delta,\epsilon} \rightsquigarrow Y_{\delta,\epsilon}$ as
$n \to\infty$. Further assume that there exists a random variable $Z$
such that for all $\delta, \epsilon> 0$, $\Pr(|Z - Y_{\delta
,\epsilon
}| > \delta) < \epsilon$. Then $Z_n \rightsquigarrow Z$ as $n \to
\infty$.
\end{lem}

The convergence of Spearman's rho is presented first.

\begin{pf*}{Proof of Proposition \ref{prop52}} Because the
complement of $\mathcal{O}$ in $[0,1]^d$ has Lebesgue measure zero, it
suffices to show that
\[
Z_n = \int_{\mathcal{O} } \mathbb{\widehat C}
^{\m}_n \,\d\Pi\rightsquigarrow Z = \int
_{\mathcal{O} } \mathbb{\widehat C} ^{\m}\,\d\Pi.
\]
Given $\delta$, $\epsilon> 0$, call on Remark \ref{rem41} to
pick $M > 0$ such that $\Pr(\|\mathbb{\widehat C}
^{\m}\|_{\mathcal{O}}
> M) <\epsilon$ and $\Pr(\|\mathbb{\widehat C} ^{\m}_n\| >M)
<\epsilon$ for all $n \in\mathbb{N}$. Then choose a compact set $K
= K_{\delta, \epsilon} \subset\mathcal{O}$ such that $\Pi
(\mathcal{O}\setminus K)< \delta/M$. Now define
\[
Y_{n,\delta,\epsilon} = \int_{K} \mathbb{\widehat
C}_n^{\m}\,\d\Pi, \qquad Y_{\delta, \epsilon} = \int
_{K} \mathbb{\widehat C}^{\m}\,\d\Pi.
\]
Theorem \ref{thm31} implies that $Y_{n,\delta,\epsilon}
\rightsquigarrow Y_{\delta, \epsilon}$ as $n \to\infty$. Furthermore,
\[
| Z_n - Y_{n,\delta,\epsilon} | = \biggl\llvert\int
_{\mathcal{O} \setminus K} \mathbb{\widehat C}_n^{\m}\,\d\Pi
\biggr\rrvert\le\bigl\|\mathbb{\widehat C}_n^{\m} \bigr\| \Pi(
\mathcal{O} \setminus K) < \frac{\delta}{M} \bigl\|\mathbb{\widehat
C}_n^{\m}\bigr\|,
\]
while
\[
| Z - Y_{\delta,\epsilon} | = \biggl\llvert\int_{\mathcal{O}
\setminus K}
\mathbb{\widehat C}^{\m}\,\d\Pi\biggr\rrvert\le\bigl\|\mathbb
{\widehat
C}^{\m}\bigr\|_{\mathcal{O}} \Pi(\mathcal{O} \setminus K) <
\frac{\delta}{M}\bigl \|\mathbb{\widehat C}^{\m}\bigr\|_{\mathcal{O}} .
\]
For all $n \in\mathbb{N}$, one then has $\Pr( \llvert Z_n -
Y_{n,\delta,\epsilon}\rrvert > \delta) \le\Pr
(\|\mathbb{\widehat C}_n^{\m}\| {\delta}/{M} > \delta) <
\epsilon$ and similarly $\Pr( \llvert Z - Y_{\delta,\epsilon}
\rrvert > \delta) < \epsilon$. The conclusion is then a
consequence of Lemma \ref{lemC1}.
\end{pf*}

\com{The following lemma, needed for the proof of Proposition \ref
{prop51}, is excerpted from Genest, Ne{\v{s}}lehov\'a and
R\'emillard \cite{GenestNeslehovaRemillard2013}.}

%
%leA.2 #&#
\begin{lem}\label{lemC2} Let $H$ be a distribution function on
$\mathbb{R}^d$ and denote by $H_n$ its empirical counterpart
corresponding to a random sample of size $n$. If the sequence of
processes $\mathbb{G}_n$ is tight with respect to the uniform norm
on the space $\mathcal{C}_b(\mathbb{R}^d )$ of bounded and
continuous functions on $\mathbb{R}^d$, then, as $n \to\infty$,
$R_n = \int\mathbb{G}_n \,\d H_n - \int\mathbb{G}_n \,\d H
\stackrel{\mathrm{ p}}{\to} 0$.
\end{lem}

\begin{pf*}{Proof of Proposition \ref{prop51}} Observe that
\[
\sqrt{n} (\tau_n-\tau) = 4\int_{\mathcal{O}} \mathbb{
\widehat C}_n^{\m}(u,v) \,\d\widehat C_n^{\m}
(u,v) + 4\int_{\mathcal{O}} \mathbb{\widehat C}_n^{\m}(u,v)
\,\d C^{\m}(u,v).
\]
First, it will be shown that, as $n\to\infty$,
\[
Z_n = \int_{\mathcal{O}} \mathbb{\widehat
C}_n^{\m}(u,v) \,\d C^{\m}(u,v)\rightsquigarrow
Z= \int_{\mathcal{O}} \mathbb{\widehat C}^{\m} (u,v) \,\d
C^{\m}(u,v).
\]
To see this, fix arbitrary $\delta,\epsilon> 0$ and use Remark
\ref{rem41} to pick $M>0$ such that $\Pr(\| \mathbb{\widehat
C}^{\m}\|_{\mathcal{O}}
> M)<\epsilon$ and $\Pr(\| \mathbb{\widehat C}_n^{\m}\| >
M)<\epsilon$
for all $n\in\mathbb{N}$. Then choose a compact set $K = K_{\delta,
\epsilon} \subset\mathcal{O}$ such that $C^{\m}(\mathcal
{O}\setminus
K) < \delta/M$. Setting
\begin{eqnarray*}
Y_{n,\delta,\epsilon} &=& \int_{K} \mathbb{\widehat
C}_n^{\m}(u,v) \,\d C^{\m}(u,v),\\
Y_{\delta,\epsilon} &=& \int_{K} \mathbb{\widehat
C}^{\m} (u,v) \,\d C^{\m}(u,v),
\end{eqnarray*}
one can invoke Theorem \ref{thm31} to deduce that
$Y_{n,\delta,\epsilon} \rightsquigarrow Y_{\delta,\epsilon}$ as
$n\to\infty$. The rest of the argument rests on Lemma \ref{lemC1},
in analogy to the proof of Proposition~\ref{prop52}.

Secondly, to establish that, as $n\to\infty$,
%
%
%eA.1 #&#
\begin{equation}
\label{eqC1} \int_{\mathcal{O}} \mathbb{\widehat C}_n^{\m}(u,v)
\,\d\widehat C_n^{\m}(u,v)\rightsquigarrow\int
_{\mathcal{O}} \mathbb{\widehat C}^{\m}(u,v) \,\d
C^{\m}(u,v),
\end{equation}
use a change of variables and the definition of $H_n^{\m}$ to write
\begin{eqnarray*}
\int_{[0,1]^2} \mathbb{\widehat C}_n^{\m}(u,v)
\,\d\widehat C_n^{\m}(u,v) &=& \int_{\mathbb{R}^2}
\mathbb{\widehat C}_n^{\m} \bigl\{F_{n1}^{\m}(x_1),F_{n2}^{\m}
(x_2) \bigr\} \,\d H_n^{\m}(x_1,x_2)
\\
&=& \int_{\mathbb{R}^2} \mathbb{G}_n(x_1,x_2)
\,\d H_n (x_1,x_2),
\end{eqnarray*}
where, for all $x_1, x_2 \in\mathbb{R}$,
\[
\mathbb{G}_n(x_1,x_2) = \int
_{[0,1]^2} \mathbb{\widehat C}_n^{\m} \bigl
\{F_{n1}^{\m}(x_1+u-1),F_{n2}^{\m}(x_2+v-1)
\bigr\} \,\d v \,\d u.
\]
It is clear that $\| \mathbb{G}_n\| \le\| \mathbb{\widehat C}_n^{\m}
\| $ and hence, by virtue of Remark \ref{rem41}, the sequence of processes
$\mathbb{G}_n$ is tight on $\mathcal{C}_b(\mathbb{R}^2)$. Lemma~\ref{lemC2}
thus implies that, as $n \to\infty$,
\[
\biggl\llvert\int_{\mathbb{R}^2} \mathbb{G}_n(x_1,x_2)
\,\d H_n (x_1,x_2) - \int
_{\mathbb{R}^2} \mathbb{G}_n(x_1,x_2)
\,\d H (x_1,x_2) \biggr\rrvert\stackrel{\mathrm{p}} {\to} 0.
\]
Undoing the change of variables and using the definitions of $H^{\m}$ and
$C^{\m}$, one finds
\begin{eqnarray*}
\int_{\mathbb{R}^2} \mathbb{G}_n(x_1,x_2)
\,\d H (x_1,x_2)& =& \int_{\mathbb{R}^2}
\mathbb{ \widehat C}_n^{\m} \bigl\{F_{n1}^{\m}(x_1),F_{n2}^{\m}
(x_2) \bigr\} \,\d H^{\m}(x_1,x_2)
\\
&=& \int_{[0,1]^2} \mathbb{\widehat C}_n^{\m}
\bigl\{F_{n1}^{\m}\circ F_1^{\m
-1}(u),F_{n2}^{\m}
\circ F_2^{\m-1}(v) \bigr\} \,\d C^{\m}(u,v)
\\
&=& \int_{\mathcal{O}} \mathbb{\widehat C}_n^{\m}
\bigl\{F_{n1}^{\m}\circ F_1^{\m
-1}(u),F_{n2}^{\m}
\circ F_2^{\m-1}(v) \bigr\} \,\d C^{\m}(u,v).
\end{eqnarray*}
Claim \eqref{eqC1} is established if one can show that, as $n \to
\infty$,
%
%
%eA.2 #&#
\begin{equation}
\label{eqC2} \bigl\| \mathbb{\widehat C}_n^{\m}
\bigl(F_{n1}^{\m}\circ F_1^{\m-1},F_{n2}^{\m}
\circ F_2^{\m-1} \bigr) - \mathbb{\widehat
C}_n^{\m}\bigr\|_K \stackrel{\mathrm{p}} {\to} 0
\end{equation}
for any fixed compact set $K \subset\mathcal{O}$. Given such a set,
one can proceed exactly as in the proof of Proposition \ref{prop44}
to find $\delta^* >0$ and a compact set $K^* \subset\mathcal{O}$ such
that for all $\delta\in(0,\delta^*)$, $K \subset K_\delta\subset K^*$.

Next, fix $\epsilon> 0$ and $\delta\in(0,\delta^*)$ and recall that
$\|\mathbb{C}^{\m}_{nj}\| \le\|\mathbb{B}_n\|$ for $j=1,2$. As in the
proof of Proposition \ref{prop44}, one has that when $\{\|\mathbb
{B}_n\| / \sqrt{n} < \delta/2\}$ holds,
\[
\bigl(F_{n1}^{\m}\circ F_1^{\m-1}(u),F_{n2}^{\m}
\circ F_2^{\m-1}(v) \bigr) = \biggl(u+ \frac{\mathbb{C}_{n1}^{\m
}(u)}{\sqrt{n}}, v+
\frac{\mathbb
{C}_{n2}^{\m}(v)}{\sqrt{n}} \biggr) \in K_\delta
\]
whenever $(u,v) \in K$. Therefore,
\[
P^* \bigl\{\bigl\| \mathbb{\widehat C}_n^{\m}
\bigl(F_{n1}^{\m}\circ F_1^{\m-1},F_{n2}^{\m}
\circ F_2^{\m-1} \bigr) - \mathbb{\widehat
C}_n^{\m}\bigr\|_K > \epsilon\bigr\}
\]
is bounded above by
\[
P^* \Bigl\{ \mathop{\sup_{(u,v),(u^*,v^*) \in K^* }}_{ |u-u^*| + |v-v^*|
< \delta}\bigl|\mathbb{
\widehat C}_n^{\m} \bigl(u^*,v^* \bigr) - \mathbb{\widehat
C}_n^{\m} (u,v)\bigr| > \epsilon\Bigr\} + P^* \bigl( \|
\mathbb{B}_n\| / \sqrt{n} \ge\delta/2 \bigr).
\]
Given that $\mathbb{\widehat C}_n^{\m}$ converges to $\mathbb
{\widehat
C}^{\m}$ on $\mathcal{C}(K^*)$,
\[
\lim_{\delta\downarrow0} \limsup_{n\to\infty} P^* \Bigl\{
\mathop{\sup_{{(u,v),(u^*,v^*) \in K^*}}}_{ |u-u^*| + |v-v^*| <
\delta} \bigl|\mathbb{\widehat
C}_n^{\m} \bigl(u^*,v^* \bigr) - \mathbb{\widehat
C}_n^{\m}(u,v)\bigr| > \epsilon\Bigr\} = 0.
\]
Claim \eqref{eqC2} now readily follows from the fact that $\|\mathbb
{B}_n\| / \sqrt{n} \stackrel{\mathrm{p}}{\to} 0$, as $n\to\infty$. In
conclusion, $\sqrt{n} (\tau_n - \tau) \rightsquigarrow\mathcal{T}_2
= 8\int_{\mathcal{O}} \mathbb{\widehat C}^{\m}(u,v) \,\d C^{\m
}(u,v)$ as $n
\to\infty$, as claimed.
\end{pf*}
\end{appendix}

\section*{Acknowledgments}
The authors thank the Editors and referees for many helpful comments.
They acknowledge grants from the Canada Research Chairs Program, the
Natural Sciences and Engineering Research Council and the Fonds de recherche du Qu\'ebec -- Nature et technologies.

%%%%%%%%%%%%%%%%%%%%%%%%%%%%%%%%%%%%%%%%%%%%%%%%%%%%%%%%%%%%%%%%%
%
% imsref loaded by akundreckaite, 2013-07-03 12:22:05
% imsref loaded by akundreckaite, 2013-07-03 12:36:09
%

%


%suskaldyti doi

\printhistory

\end{document}